\documentclass[11pt,a4paper]{article}
\usepackage{amsmath,amsthm,amssymb}
\usepackage{nicefrac}
\usepackage{amssymb}
\usepackage{tikz,pgfplots}
\usetikzlibrary{plotmarks,external,spy}
\usepgfplotslibrary{colormaps}
\pgfplotsset{compat=1.10}
\usepackage{subcaption}
\usepackage{hyperref}
\usepackage{placeins}
\usepackage[utf8]{inputenc}
\usepackage{graphicx}
\usepackage{siunitx}
\definecolor{winered}{rgb}{0.5,0,0}
\definecolor{darkgreen}{rgb}{0.0,0.5,0}
\hypersetup{
colorlinks=true,
linkcolor=winered,
citecolor=darkgreen
}
\usepackage{titling}
\usepackage{mathrsfs}
\usetikzlibrary{plotmarks}
\pgfplotsset{colormap={ngsolve}{rgb255=(0,0,255) rgb255=(0,255,255) rgb255=(0,255,0) rgb255=(255,255,0) rgb255=(255,0,0)}}
\newcommand{\xvec}{\mathbf{x}}

\newcommand{\vvec}{\mathbf{v}}
\newcommand{\vbar}{\bar{\mathbf{v}}}
\newcommand{\vhat}{\hat{\mathbf{v}}}
\newcommand{\Vvec}{\mathbf{V}}
\newcommand{\wvec}{\mathbf{w}}
\newcommand{\evec}{\mathbf{e}}

\newcommand{\nvec}{\mathbf{n}}

\newcommand{\qvec}{\mathbf{q}}

\newcommand{\ip}{\text{ip}}

\newcommand{\R}[1]{\mathbb{R}^{#1}}

\newcommand{\rhozo}{\rho_{\Vvec,T}}
\newcommand{\Vzo}{\Vvec_{\Vvec,T}}
\newcommand{\Tzo}{T_{\Vvec,T}}

\newcommand{\fzo}{f^{\Vvec,T}}
\newcommand{\fzon}{f^{\Vvec_n,T_n}}
\newcommand{\phizo}{\varphi^{\Vvec,T}}

\newcommand{\diver}{\text{div}}
\newcommand{\refer}[1]{(\ref{#1})}
\newcommand{\ndof}{\text{ndof}}
\newcommand{\ndofv}{\text{ndof}_\vvec}

\newcommand{\el}{\mathcal{T}}

\newcommand{\mesh}{\mathscr{T}}
\newcommand{\intrd}{\int\limits_{\R{3}}}

\newcommand{\fzoup}{f^{\Vvec,T,\text{up}}}
\newcommand\restr[2]{{
  \left.\kern-\nulldelimiterspace 
  #1 
  \vphantom{\big|} 
  \right|_{#2} 
  }}

\DeclareMathOperator{\tr}{tr}
\newtheorem{theorem}{Theorem}

\newtheorem{proposition}[theorem]{Proposition}

\newtheorem{remark}[theorem]{Remark}

\title{A spatial discontinuous Galerkin method with rescaled velocities for the Boltzmann equation}
\author{Gerhard Kitzler and Joachim Sch\"oberl}

\begin{document}
\maketitle
\begin{abstract}
\noindent In this paper we present a numerical method for the Boltzmann equation. It is a spectral discretization in the velocity and a discontinuous Galerkin discretization in physical space. To obtain uniform approximation properties in the mach number, we shift the velocity by the (smoothed) bulk velocity and scale it by the (smoothed) temperature, both extracted from the density distribution. The velocity trial functions are polynomials multiplied by a Maxwellian. Consequently, an expansion with a low number of trial functions already yields satisfying approximation quality for nearly equilibrated solutions. By the polynomial test space, density, velocity and energy are conserved on the discrete level. Different from moment methods, we stabilize the free flow operator in phase space with upwind fluxes. Several numerical results are presented to justify our approach.
\end{abstract}

\vspace{-.5cm}
\section{The Boltzmann equation}
To begin with, let us consider a rarefied gas. The rarefaction shall be such that -- due to the lack of sufficient particle interactions -- the gas is locally not in thermal equilibrium. On the other hand, collisions between particles are still important to describe the flow accurately. A mathematical description for such a gas is given by the Boltzmann transport equation, which provides a statistical description. In addition to rarefied gas dynamics, kinetic models as the Boltzmann equation have among others, attracted a wide range of applications including plasma physics, electron transport in semiconductors and also disciplines from biology. Typically, the unknown in kinetic equations is denoted as $f = f(t,\xvec,\vvec), \xvec \in \Omega, \vvec \in \R3$. $\Omega\subset \R3$ is the physical domain under investigation, $f$ is the density in phase, in other words
\begin{equation}
\rho(t,\xvec) := \intrd f(t,\xvec,\vvec)\,d\vvec
\label{eq:mac_rho}
\end{equation}
gives the physically measurable mass density $\rho(t,\xvec)$ at time $t$ and position $\xvec$. Similar, the mean velocity $\Vvec(t,\xvec)$, the energy density $E(t,\xvec)$, the stress tensor $P$, the (scalar) pressure $p(t,\xvec)$ and the temperature $T(t,\xvec)$ are derived from moments of $f$:
\begin{equation}
\begin{alignedat}{2}
\mathbf{V}(t,\xvec)    &:= \frac{1}{\rho(t,\xvec)}\intrd \vvec f(t,\xvec,\vvec)\,d\vvec & \quad \quad P(t,\xvec) &:= \intrd (\vvec-\Vvec)\otimes (\vvec-\Vvec)f(t,\xvec,\vvec)\,d\vvec \\
E(t,\xvec) &:= \frac{1}{2}\intrd |\vvec|^2f(t,\xvec,\vvec)\,d\vvec & \quad \quad p(t,\xvec) &:= \frac{1}{3}\tr{P}(t,\xvec) \\  
\qvec(t,\xvec) &:= \frac{1}{2}\intrd |\vvec|^2 \vvec f(t,\xvec,\vvec)\,d\vvec & \quad\quad T(t,\xvec)&:= \frac{1}{\rho(t,\xvec)}p(t,\xvec).
\end{alignedat}
\label{eq:mac_remainders}
\end{equation}

In contrast to $f$, these quantities are measurable in experiments. \par
The time evolution of $f$ is governed by the Boltzmann equation:
\begin{equation}
\frac{\partial }{\partial t}f + \diver_\xvec(\vvec f) = Q(f), \quad \xvec\in \Omega,\;\vvec \in \R3
\label{eq:boltzmann_equation}
\end{equation}
The left hand side of this equation takes into account free flow of the molecules. The right hand side models binary collisions, consequently, $Q$ is termed a collision operator. In case of the Boltzmann equation it has the specific form
\begin{equation}
Q(f(t,\xvec,\,.\,) )(\vvec) := \intrd\int\limits_{S^2}B(\vvec,\wvec,\evec') [f(t,\xvec,\vvec')f(t,\xvec,\wvec')-f(t,\xvec,\vvec)f(t,\xvec,\wvec)]d\evec'd\wvec,
\label{eq:collop}
\end{equation}
with the pre collision velocities given by 
\begin{equation}
\vvec' = \frac{\vvec+\wvec}{2}+\evec'\frac{|\vvec-\wvec|}{2}\;\text{ and }\; \wvec' = \frac{\vvec+\wvec}{2}-\evec'\frac{|\vvec-\wvec|}{2}.
\label{eq:precollvel}
\end{equation}
The function $B(\vvec,\wvec,\evec') = B(|\vvec-\wvec|,\tfrac{(\vvec-\wvec)\cdot e'}{|\vvec-\wvec|})$ is called a collision kernel. It is the probability density for two particles with velocity $\vvec'$ and $\wvec'$ to collide and result in velocities $\vvec$ and $\wvec$. In many relevant cases, $B$ factorizes as 
\begin{equation}
B(\vvec,\wvec,\evec') = b_r(|\vvec-\wvec|)b_\theta(\tfrac{(\vvec-\wvec)\cdot \evec'}{|\vvec-\wvec|}), \quad \quad b_r(|\vvec-\wvec|) = |\vvec-\wvec|^\beta,\quad\beta \in [0,1].
\label{eq:collisionkernel}
\end{equation}
The collision operator $Q$ is non linear, local in time and position, but global in velocity, i.e. $Q(f)(t,\xvec,\vvec) = Q(f(t,\xvec,\,.\,) )(\vvec)$. It satisfies conservation of mass, momentum and energy, what is expressed by the equations
\begin{equation}
\intrd Q(f)\left(\begin{matrix}1\\ \vvec \\ |\vvec|^2 \end{matrix}\right) = \mathbf{0}.
\label{eq:coll_invariants}
\end{equation}
The functions $1,\vvec,|\vvec|^2$ are consequently termed collision invariants. Also, each $\phi$ satisfying $\int Q(f)\phi\,d\vvec = 0$ is a linear combination of the basic collision invariants.
\par
Concerning our numerical method also the null space $\ker Q$, which is spanned by the Maxwell distributions, is of importance:
\begin{equation}
\mathrm{ker}\,Q = \left\{f: f(t,\xvec,\vvec) = \rho(t,\xvec) e^{-\left|\frac{\vvec-\Vvec(t,\xvec)}{\sqrt{T(t,\xvec)}}\right|^2} \right\}
\label{eq:collisionnullspace}
\end{equation}
The quantities $\rho, T$ and $\Vvec$ in \refer{eq:collisionnullspace} are in accordance with the definition of mass, bulk velocity, and temperature in \refer{eq:mac_rho} and \refer{eq:mac_remainders}. Clearly, each non negative stationary state $f_\infty$ for the space homogeneous problem $\partial_t f = Q(f)$ is in $\ker Q$.

The conservation properties \refer{eq:coll_invariants} are also at the basis of the macroscopic behaviour: Testing \refer{eq:boltzmann_equation} with the collision invariants yields an (unclosed) system of equations for the macrscopic quantities defined in \refer{eq:mac_rho} and \refer{eq:mac_remainders}:
\begin{equation}
\begin{alignedat}{2}
&\frac{\partial}{\partial t}\rho + \diver(\rho \Vvec) & &= 0 \\
&\frac{\partial}{\partial t}(\rho\Vvec) + \diver(P + \rho \Vvec \otimes \Vvec) & &= 0 \\
&\frac{\partial}{\partial t}(3p + |\Vvec|^2\rho) + \diver(2\qvec + 3p+|\Vvec|^2 \rho \Vvec + P \Vvec) & &= 0.
\end{alignedat}
\label{eq:unclosed_cons_law}
\end{equation}
\par
Finally, we discuss boundary conditions for the Boltzmann equation. 
Let us fix notation: We denote the normal vector at a point $\xvec \in \partial \Omega$ by $\nvec$. Given the vector $\nvec$, we define $\R3_\text{in} := \{\vvec \in \R3: \vvec\cdot \nvec<0\}$ and $\R3_\text{out} := \R3\setminus \R3_\text{in}$ as the in and outgoing directions respectively. Note that $\nvec,\,\R3_\text{in}$ as well as $\R3_\text{out}$ depend on $\xvec$.
\begin{itemize}
\item Inflow
\begin{equation}
f(t,\xvec,\vvec) = f_\text{in}(t,\xvec,\vvec),\quad\quad \xvec\in (\partial \Omega)_{\text{in}},\vvec \in \R3_\text{in}(\xvec),
\label{eq:bc_inflow}
\end{equation}
i.e. at a fixed $\xvec \in \partial \Omega$, $f$ is prescribed at  velocities $\vvec \in \R3$, which point to the inside of $\Omega$.
\item Specular reflection
\begin{equation}
f(t,\xvec,\vvec) = f(t,\xvec,\vvec-2(\vvec\cdot \nvec)\nvec),\quad\quad \xvec (\partial \Omega)_{\text{spec}},\vvec \in \R3_\text{in}(\xvec),
\label{eq:bc_specular}
\end{equation}
i.e. a molecule hitting the boundary behaves like a billard ball.
\item Diffuse reflection
\begin{equation}
f(t,\xvec,\vvec) = c e^{-\left|\frac{v-\Vvec_\text{bnd}}{\sqrt{2T_\text{bnd}}} \right|^2},\quad\quad \xvec\in (\partial \Omega)_{\text{diff}},\vvec \in \R3_\text{in}(\xvec),
\label{eq:bc_diffuse}
\end{equation}
i.e. a particle hitting the boundary is redistributed from a Maxwell distribution with bulk velocity $\Vvec_\text{bnd}$ and temperature $T_\text{bnd}$.
The constant $c$ is a normalization to guarantee vanishing normal flux across the boundary, i.e.
$$c = \Biggl(\;\,\int\limits_{\R3_\text{out}(\xvec)}\!\!\!f(t,\xvec,\vvec)\vvec\cdot \nvec\,d\vvec \Biggr) \Biggl(\;\,\int\limits_{\R3_\text{in}(\xvec)} e^{-\left|\frac{v-\Vvec_\text{bnd}}{\sqrt{2T_\text{bnd}}} \right|^2}| \vvec\cdot \nvec | \,d\vvec\Biggr)^{-1}.$$
\end{itemize}
For a detailed discussion about the Boltzmann equation we refer the reader to \cite{Cercignani_1994}.
\par
Now that we turn to numerical methods, let us summarize some of the challenges that need to be addressed:
\begin{itemize}
\item \textbf{high dimensionality}
The method uses a tensor product in $\xvec$ and $\vvec$. Additionally, the three dimensional velocity domain $\R3$ provides -- due to its Cartesian structure -- the possibility for additional factorization. This does of course not reduce the number of degrees of freedom, but it reduces the number of flops required for function evaluation, integration, etc.
\item \textbf{non linearity and high compuational costs for $Q$} We refer to the space homogeneous paper \cite{Kitzler2018} where a reduction of computational costs for the collision operator based on efficient basis transformations is presented.
\item \textbf{unbounded domain} The method uses global basis functions w.r.t $\vvec$. The arising integrals in the variational formulation exist and can be evaluated exactly using Gauss Hermite quadrature formulas.
\item \textbf{stabilization of the advection operator} We apply standard upwinding to the advection operator to improve stability.
\item \textbf{Large variations of $f(t,\xvec,.)$'s support} This requires either a large number of degrees of freedom w.r.t. the velocity variable or adaptivity w.r.t. $f$'s support. We present an approach based on resacling and shifting $f$'s velocity variable.
\end{itemize}
According to the high computational costs for the collision operator, stochastic methods as presented in \cite{doi:10.1143/JPSJ.49.2050,MR845926,MR1352466,Rjasanow_2005} are of great interest since their computational costs are linear w.r.t the number of degrees of freedom. In contrast to deterministic methods, the solutions are subject to stochastic noise, the convergence in general is rather slow. However, discontinuouities w.r.t. the velocity are resolved very well. On the other hand, investigating deterministic methods, a considerable portion of methods is based on trigonometric spectral techniques, see \cite{Pareschi_fft_homo,MR1439069,Bobylev_hard_sphere,MR1756425,GAMBA20092012,Ibragimov2002,
10.2307/43693916,MR2746671,doi:10.1080/00411450008205876}. The density distribution is expanded to a trigonometric polynomial on some bounded domain. The orthogonality of the trigonometric expansion then gives efficient algorithms for the collision operator. By further approximation of the resulting collision weights a convolution structure is obtained, leading to the classical fast spectral methods, see exemplary \cite{MR2240637,filbetnlogn,doi:10.1137/16M1096001}. The drawback of these methods is that the integrals defining the collision operator, as well as the approximation domain have to be truncated. The periodicity of the trigonometric basis combined with the larger support of $Q$ regarding the density distribution $f$ may additionally produce aliasing errors. In fourier spectral methods,  typically mass is conserved, momentum and energy up to spectral accuracy. This is a consequence of missing collision invariants in the test space.
In recent years the classical Fourier spectral methods were added polynomial spectral methods. In \cite{GambaR2017,Torsten2017} I. Gamba, S. Rjasanow and T. Ke\ss{}ler present a similar approximation in the velocity variable combined with a finite volume method for the spatial transport. Also in the group of Hiptmair a similar discretization was developed, to the best of our knowledge, this method is currently only available in two velocity dimensions, see \cite{FGH14_563,GHP15_628,MR1701710}. 
Let us also mention discrete velocity methods. Simplified speaking, particles are only allowed to have velocities on a finite grid. Consequently, the collision mechanism needs to be discretized in such a way, that pre and post collision velocities are nodes of the grid, while maintaining the main physical properties of the collision process. We refer to \cite{bobylevdvm1995,doi:10.1080/00411459608204829,Panferov_DVM}.

\par
The remainder of this paper is outlined as follows: In section \ref{subsec:normalized_density} we introduce the standardized distribution function and show our approximation to that function in sections \ref{subsec:vspace} - \ref{subsec:xvspace}. Then, in section \ref{sec:variational_formulation} we present our discontinuous Galerkin method, with a rescaled velocity variable. Finally, section \ref{sec:numerical_results} provides numerical results as a validation for our method.

\section{Approximating the density distribution}
\subsection{The standardized distribution function}
\label{subsec:normalized_density}
In spatially inhomogeneuous problems, large variations in the gas' mean velocity $\Vvec_f$ and temperature $T_f$ may occur. That requires -- when the same approximation w.r.t. the velocity is used at each spatial point -- a larger number of velocity degrees of freedom and consequently a higher computational effort to obtain reasonable accuracy. In order to reduce the numerical work, we propose to introduce a standardized distribution $\fzo$ and test function $\phizo$, defined by
$$
\fzo(t,\xvec,\vvec):=f(t,\xvec,\sqrt{T}\vvec+\Vvec),\quad\quad \phizo(\xvec,\vvec):=\varphi(\xvec,\sqrt{T}\vvec+\Vvec).
$$
The macroscopic properties density, velocity and temperature of $\fzo$ are denoted by $\rhozo$, $\Vzo$ and $\Tzo$ respectively, those for $f$ are denoted with a subscript $f$.
The fields $T = T(t,\xvec)$ and $\Vvec=\Vvec(t,\xvec)$, denoted as ansatz temperature and velocity respectively are here assumed to be given functions on the spatial domain $\Omega$. For good approximation properties, they should be close to the gas mean temperature $T_f(t,\xvec)$ and velocity $\Vvec_f(t,\xvec)$. Instead of solving the Boltzmann equation for $f$, we  derive an equation for $\fzo$. Contrary to $f$, this function has almost zero velocity and a  temperature close to one. As a consequence the same approximation space w.r.t. the velocity gives equal approximation properties for different $\Vvec$ and $T$. 
\par 
In the following propositions, we collect relations between $f$ and $\fzo$. Let us begin with the relation between the macroscopic properties of $f$ and $\fzo$:
\begin{proposition}
The macroscopic properties of $f$ and $\fzo$ are related via
\begin{equation*}
\begin{aligned}
\rhozo &= T^{-\nicefrac{3}{2}}\rho_f \\
\Vzo &= \frac{1}{\sqrt{T}} (\Vvec_f - \Vvec) \\
\Tzo &= \frac{T_f}{T}.
\end{aligned}
\end{equation*}
\label{prop:macroscopics_of_fzo}
\end{proposition}
The next proposition shows how the collision operator is expressed in terms of the standardized density:
\begin{proposition}
Assume the collision kernel $B = B(\vvec,\wvec,\evec)$ is of the form \refer{eq:collisionkernel}. Then there holds 
$$
\int\limits_{\R3}Q(f)\varphi \,d\vvec = T^{3+\nicefrac{\beta}{2}}\int\limits_{\R3}Q(\fzo)\phizo \,d\vvec
$$
for all suitable test functions $\varphi(\vvec)$.
\label{prop:collision_for_fzo}
\end{proposition}
\begin{proof}
For the proof we start with the left hand side. We perform the change of variables $e'\mapsto e',\,\vvec\mapsto\sqrt{T}\vvec+\Vvec,\,\sqrt{T}\wvec\mapsto \wvec+\Vvec$. The determinant is $T^{3}$, the collision kernel transforms as $B(\vvec,\wvec,\evec') \mapsto \sqrt{T}^\beta B(\vvec,\wvec,\evec')$. Using the expression for the pre collisional velocities \refer{eq:precollvel}, we note that
they transform according to $\vvec' \mapsto\sqrt{T}\vvec'+\Vvec$ and $\wvec'\mapsto \sqrt{T}\wvec'+\Vvec$ into transformed pre collision velocities.
\end{proof}
\begin{remark}
The mass integrals as well as the fluxes transform in a similar way as the collision operator. Applying the same transformation as in the previous proof we obtain
$$
\intrd f\phi \, d\vvec = T^{\nicefrac{3}{2}}\intrd \fzo \phizo\,d\vvec,\quad\quad \intrd \vvec f\phi \, d\vvec = T^{\nicefrac{3}{2}}\intrd (\sqrt{T}\vvec+\Vvec)\fzo \phizo\,d\vvec.
$$
\label{rem:transformedmassandflux}
\end{remark}
\subsection{Velocity discretization}
\label{subsec:vspace}
Let us now discuss the discretization in the velocity domain: We consider a function $\fzo = \fzo(t,\vvec)$ to motivate the choice of trial and test space for the velocity discretization. Since $Q$ is global in $\vvec$, a global trial space in the velocity is suitable. Close to a local equilibrium, $f$ is close to a Maxwellian with the gas' mean density, velocity and temperature \refer{eq:collisionnullspace}. Consequently, we expect excellent approximation properties if we approximate the density with polynomials, multiplied by an appropriate Maxwellian. This approximation also ensures that the kernel $\ker Q$ and consequently the local thermal equilibria are in the trial space. Let us fix notation: The trial space is
\begin{equation*}
V_N:=\left\{f \in L_2(\R3): f(\vvec) = e^{-|\vvec|^2}p(\vvec):\,p \in Q^N(\R3) \right\},
\end{equation*}
where $Q^N(\R3)$ is the space of polynomials of partial degree $N$ on $\R3$. \par
The test space is motivated by the conservation properties of $Q$, \refer{eq:coll_invariants} which lead to the macroscopics equations \refer{eq:unclosed_cons_law}. These local conservations result from testing $Q$ with the collision invariants $1, \vvec, |\vvec|^2$. In other words, they are naturally satisfied in a Galerkin projection if the collision invariants are in the test space. Accordingly, the test space is defined as
\begin{equation*}
W_N:=Q^N(\R3).
\end{equation*}
Note that $N\geq 2$ is a minimum requirement to satisfy \refer{eq:coll_invariants}.
\begin{remark}
In the simplified situation of the space homogeneuous problem $\partial_t f = Q(f)$, the choice of the ansatz parameters $\Vvec$ and $T$ is straight forward. Since $\Vvec_f$ and $T_f$ are conserved over time, we calculate them a priori from the initial density and set $\Vvec = \Vvec_f$ as well as $T = T_f$. This is described in detail in \cite{Kitzler2018}, were we  studied the approximation properties of the spaces $V_N$ and $W_N$ in space homogeneuous problems.
\end{remark}
\textbf{Polynomial basis in $Q^N(\R3)$.}
Let us now present the polynomial basis of $Q^N(\R3)$ we use. $Q^N(\R3)$ is the space of polynomials of partial degree $N$. To define a suitable basis of $Q^N(\R3)$ we denote by $(v_\text{ip},\omega_\text{ip}),\,\text{ip}=0\ldots N$ the nodes and weights of a Gauss-Hermite quadrature \cite{shen2011spectral}, satisfying
\begin{equation*}
\int\limits_{\R{}} e^{-v^2}p(v) = \sum\limits_{\ip=0}^{N} \omega_\ip p(v_\ip),\;\forall\, p \in P^{2N+1}(\R{}).
\end{equation*}
\par
In one dimension, we define Lagrange collocation polynomials to the Gauss-Hermite quadrature nodes and denote them by $l$
\begin{equation*}
l_j(v) := \prod\limits_{\genfrac{}{}{0pt}{1}{i=0}{i\neq j}}^{N} \frac{v-v_i}{v_j-v_i}.
\end{equation*}
The three dimensional basis is obtained as the tensor product of the 1D polynomials and its elements are denoted by $L_j,j=0\dots \ndofv-1=(N+1)^3-1$. Here, $\ndof_\vvec$ denotes the number of degrees of freedom.
\begin{equation*}
L_j(\vvec) := l_l(\vvec_1)l_m(\vvec_2)l_n(\vvec_3),\quad j = (N+1)^2l + (N+1)m+n.
\end{equation*}
The Gauss-Hermite quadrature is exact for the integrals in \refer{eq:orthogonality_basis_lagrange}. Thus, the Lagrange polynomials provide an orthogonal basis of $Q_N(\R3)$ equipped with the Maxwellian weighted $L_2$ inner product. In addition to a diagonal mass matrix, also the flux matrix is diagonal:
\begin{equation}
\begin{alignedat}{2}
M^{\vvec}_{nm}&=\intrd e^{-|\vvec|^2}L_m(\vvec)L_n(\vvec)\,d\vvec  & & = \sum\limits_{\ip=0}^{\ndofv-1}\omega^{(3)}_\ip L_m(\vvec_\ip^{(3)})L_n(\vvec_\ip^{(3)})= \delta_{m,n}\omega^{(3)}_n\\
F^{\vvec}_{nm}&=\intrd \vvec e^{-|\vvec|^2}L_m(\vvec)L_n(\vvec)\,d\vvec & & = \sum\limits_{\ip=0}^{\ndofv-1}\omega^{(3)}_\ip \vvec_\ip^{(3)} L_m(\vvec_\ip^{(3)})L_n(\vvec_\ip^{(3)}) = \delta_{m,n}\vvec_{n}^{(3)}\omega^{(3)}_n
\end{alignedat}
\label{eq:orthogonality_basis_lagrange}
\end{equation}

Note that in \refer{eq:orthogonality_basis_lagrange}, $\vvec_{\ip}^{(3)}$ and $\omega^{(3)}_{\ip}$ are nodes and weights of a Cartesian product of Gauss Hermite formulas with nodes $\vvec_{\ip}^{(3)} = (v_i,v_j,v_k)$ and $\omega^{(3)}_{\ip} = \omega_i\omega_j\omega_k$, where $\ip=(N+1)^2i+(N+1)j+k$.
\subsection{Space-velocity discretization}
\label{subsec:xvspace}
Let us consider the approximation of $\fzo$ in phase space. We approximate the density distribution in the spatial variable by piecewise, discontinuous polynomials. To that end we consider a given mesh $\mesh$ of the spatial domain $\Omega$, i.e. $\Omega = \bigcup_{\el \in \mesh}\el$. The phase space $\Omega \times \R3$ is discretized as $\bigcup\limits_{\el \in \mesh} (\el \times \R3)$. 
Using the mesh $\mesh$, we define trial and test space as tensor products of a standard Discontinuous Galerkin space
$$
V_h^{\text{DG}}:=\{u\in L_2(\Omega): \restr{u}{\el} \in P^N(\el),\,\forall\, \el \in \mesh\}
$$
and the spaces $V_N$ and $W_N$ introduced in section \ref{subsec:vspace}. Thus, we have 
$$
V_{h,N} := V_h^{\text{DG}} \otimes V_N,\quad \quad\quad W_{h,N} := V_h^{\text{DG}} \otimes W_N.
$$
The expansion of $\fzo$ reads
$$
\restr{\fzo}{\el}(t,\xvec,\vvec) = e^{-\left|\vvec\right|^2}\sum\limits_{i,j} \phi_i(\xvec) c_{i,j}(t) L_j\left(\vvec\right),\quad \el \in \mesh,
\label{eq:expansion_of_fzo}
$$
with time dependent coefficient matrix $c_{ij}(t)$ and spatial basis functions $\phi_i, i=0\ldots \ndof_{\el_i}$, where $\ndof_{\el_i}$ denotes the number of spatial degrees of freedom on the $i-$th element. The basis functions for different elements $\el_1$ and $\el_2$ are uncoupled, consequently the resulting mass matrix is block diagonal only. Thus, to apply its inverse, only local mass matrices need to be inverted. Using a lexicographic enumeration of the degrees of freedom in $V_{h,N}$ and $W_{h,N}$, the $(\xvec,\vvec)$ mass matrix is the Kronecker product of the spatial and the velocity mass matrices $M^\xvec$ and $M^\vvec$ respectively:
$$
M^{\xvec,\vvec}=M^\xvec \otimes M^\vvec
$$
The inverse of $M^{\xvec,\vvec}$ is the Kronecker product of the inverse matrices:
$$
(M^{\xvec,\vvec})^{-1}=(M^\xvec)^{-1} \otimes (M^\vvec)^{-1}
$$
Its application to a given coefficient matrix $c \in \R{\ndof_\xvec \times \ndofv}$ can be written as $(M^\xvec)^{-1} c (M^\vvec)^{-1}$, so we first scale the $i-$th column of $c$ with the $i-$th entry of the inverse velocity mass matrix and then solve with the spatial mass matrtix.
\begin{remark}
In non curved elements the spatial mass matrix $M^{\xvec}$ is diagonal only. $M^{\xvec,\vvec}$ and its inverse are diagonal for such elements. The $j$-th row of $c$ is scaled with the $j$-th entry of the matrix $(M^\xvec)^{-1}$ and the $i$-th column with the $i$-th entry of $(M^{\vvec})^{-1}$.
\end{remark}
\section{\texorpdfstring{The spatial discontinuous Galerkin method with rescaled velocities}{DGRVM}}
\label{sec:variational_formulation}
Let us now derive a variational formulation for $\fzo$. To that end let us write the Boltzmann equation for $f$ in weak form. As usual, we multiply \refer{eq:boltzmann_equation} by an appropriate test function $\varphi = \varphi(\xvec,\vvec)$ and integrate w.r.t space and velocity. We slightly abuse notation and write $\el$ for the phase space elements $\el \times \R3$, as well as $\mesh$ for the mesh of the phase space $\Omega \times \R3$. On each element $\el\in \mesh$ we integrate the divergence by parts, resulting in a facet integral and an additional element part.
Summing w.r.t. the elements yields the following integral equation for $f$:
\begin{equation}
\sum\limits_{\el \in \mesh}\int\limits_{\el} \partial_t f\varphi\,d(\xvec,\vvec) + \int\limits_{\partial \el}\vvec \cdot \nvec f\varphi \,d(\xvec,\vvec) - \int\limits_{\el} \vvec f\cdot \nabla_\xvec \varphi\,d(\xvec,\vvec) -\int\limits_{\el} Q(f)\varphi \,d(\xvec,\vvec) = 0
\label{eq:variational_formulation}
\end{equation}
Now in order to arrive at an equation for $\fzo$, we plug its definition into \refer{eq:variational_formulation}, use proposition \ref{prop:collision_for_fzo} and remark \ref{rem:transformedmassandflux} to express the collision and facet integrals in terms of $\fzo$ and $\phizo$:
\begin{equation}
\begin{alignedat}{1}
\sum\limits_{\el \in \mesh}&\int\limits_{\el} \partial_t \fzo(t,\xvec,\tfrac{\vvec-\Vvec}{\sqrt{T}})\phizo(\xvec,\tfrac{\vvec-\Vvec}{\sqrt{T}})\,d(\xvec,\vvec) + \int\limits_{\partial \el}T^{\nicefrac{3}{2}}(\sqrt{T}\vvec+\Vvec) \cdot \nvec \fzo\phizo\,d(\xvec,\vvec) - \\ & \int\limits_{\el} \vvec \fzo(t,\xvec,\tfrac{\vvec-\Vvec}{\sqrt{T}}) \cdot \nabla_\xvec \left[\phizo(\xvec,\tfrac{\vvec-\Vvec}{\sqrt{T}})\right]\,d(\xvec,\vvec) -  \int\limits_{\el} T^{\nicefrac{\beta}{2}+3} Q(\fzo)\phizo\,d(\xvec,\vvec) = 0
\end{alignedat}
\label{eq:variational_formulation_intermediate}
\end{equation}
To improve stability, we use upwind fluxes in the facet integrals: Let us consider an arbitrary edge $e \in \partial \el, \el \in \mesh$, connecting elements $\el$ and $\tilde{\el}$. On the edge $e$, the upwind value is defined by 
\begin{equation}
\restr{\fzoup}{e}(t,\xvec,\vvec) := \begin{cases} \restr{\fzo(t,\xvec,\vvec)}{\el} & \text{if } (\sqrt{T}\vvec+\Vvec)\cdot \nvec>0\\ \restr{\fzo(t,\xvec,\vvec)}{\tilde{\el}} & \text{else} \end{cases}.
\label{eq:defupwind}
\end{equation}
Note that -- in contrast to classical moment methods -- the above upwind decision is done in phase space. For facets $e \in \partial \Omega$, the upwind value is defined as in \refer{eq:defupwind}, with the only difference that the function values from the neighbouring elements are replaced by the boundary values given by  \refer{eq:bc_inflow} - \refer{eq:bc_diffuse}. \par
Next we calculate the $\xvec$ gradient of the test function $\phizo$. 
By the chain rule there holds
\begin{equation*}
\nabla_\xvec\left[ \phizo(\xvec,\tfrac{\vvec-\Vvec}{\sqrt{T}})\right] = (\nabla_\xvec \phizo)(\xvec,\tfrac{\vvec-\Vvec}{\sqrt{T}}) - \frac{1}{T}\left(  \tfrac{1}{2} \nabla_\xvec T \otimes \tfrac{\vvec-\Vvec}{\sqrt{T}}  + \nabla_\xvec \Vvec\sqrt{T}\right)(\nabla_\vvec \varphi)(\xvec,\tfrac{\vvec-\Vvec}{\sqrt{T}}),
\label{eq:xgradientphi}
\end{equation*}
with $(\nabla_\xvec \Vvec)_{ij} = \partial_{\xvec_i} \Vvec_j$. Consequently, if we use the expression for the $\xvec$ gradient of the test function and perform the subsitution $\tfrac{\vvec-\Vvec}{\sqrt{T}} = \tilde{\vvec}$ we can replace the term involving the $\xvec$ gradient in \refer{eq:variational_formulation_intermediate} by
\begin{equation*}
\int\limits_{\el} \vvec \fzo(t,\xvec,\tfrac{\vvec-\Vvec}{\sqrt{T}}) \cdot\nabla_\xvec\left[ \phizo(\xvec,\tfrac{\vvec-\Vvec}{\sqrt{T}})\right]\,d(\xvec,\vvec) = \int\limits_{\el} T^{\nicefrac{3}{2}}\fzo \mathbb{D}_\xvec\nabla \phizo\,d(\xvec,\vvec).
\end{equation*}
The coefficient $\mathbb{D}_\xvec \in \R{1\times 6}$ is $\mathbb{D}_\xvec:=(\sqrt{T}\vvec+\Vvec)^T\left(\begin{matrix}\text{Id}_3 & - \frac{1}{T}\left(  \tfrac{1}{2} \nabla_\xvec T \otimes \vvec  + \nabla_\xvec \Vvec\sqrt{T}\right) \end{matrix}\right)$.
\par
Next we consider the time derivative of $\fzo$ in \refer{eq:variational_formulation_intermediate}. Let us integrate \refer{eq:variational_formulation_intermediate} w.r.t. time over $[t_0,t_{n_t}]$ and apply integration by parts to the time derivative. For a fixed element $\el \in \mesh$ this gives
\begin{equation*}
\begin{aligned}
\int\limits_{t_0}^{t_{n_t}}\int\limits_\el \partial_t \fzo(t,\xvec,\tfrac{\vvec-\Vvec}{\sqrt{T}})&\phizo(\xvec,\tfrac{\vvec-\Vvec}{\sqrt{T}}) d(\xvec,\vvec)dt = \\ \int\limits_\el T^{\nicefrac{3}{2}}\fzo &\bigg.\phizo\,d(\xvec,\vvec) \bigg|_{t=t_0}^{t=t_{n_t}}  - \int\limits_{t_0}^{t_{n_t}}\int\limits_\el T^{\nicefrac{3}{2}}\fzo \mathbb{D}_t\nabla_\vvec \phizo \,d(\xvec,\vvec)dt.
\end{aligned}
\end{equation*}
The $1\times 3$ coefficient vector $\mathbb{D}_t$ is given as $\mathbb{D}_t:=\left(\frac{\partial_t V}{\sqrt{T}} + \frac{1}{2T}\vvec\partial_t T\right)^T$. 
Now let us finally combine $\mathbb{D}_\xvec$ and $\mathbb{D}_t$ to write the equation in compact form, we let $\mathbb{D}:=\mathbb{D}_\xvec + (0,0,0,\mathbb{D}_t)$. Plugging everything together we obtain an equation for the standardized density $\fzo$:
\begin{equation}
\begin{aligned}
&\text{Find } \fzo \in L_2([t_0,t_{n_t}],V_{h,N}) \text{ s.t. }\forall \phizo \in W_{h,N}: \\
&\sum\limits_{\el \in \mesh}\left.\int\limits_\el T^{\nicefrac{3}{2}}\fzo\phizo d(\xvec,\vvec)\right|_{t=t_0}^{t=t_{n_t}} + \int\limits_{t_0}^{t_{n_t}}\int\limits_{\partial \el}T^{\nicefrac{3}{2}}(\sqrt{T}\vvec+\Vvec) \cdot \nvec \fzoup\phizo\,d(\xvec,\vvec)dt - \\ & \int\limits_{t_0}^{t_{n_t}}\int\limits_{\el} T^{\nicefrac{3}{2}}\fzo \mathbb{D}\nabla \phizo\,d(\xvec,\vvec)dt -  \int\limits_{t_0}^{t_{n_t}}\int\limits_{\el} T^{\nicefrac{\beta}{2}+3} Q(\fzo)\phizo\,d(\xvec,\vvec)dt = 0
\end{aligned}
\label{eq:variational_formulation_zo_final}
\end{equation}
\subsection{Time discretization}
Let us next arrive at a time discretization for \refer{eq:variational_formulation_zo_final}. To that end we introduce the following notation: We consider a subdivision of $[t_0,t_{n_t}] = \bigcup_{i=0}^{n_t-1} [t_i,t_{i+1}]$ and denote by $c^{n}$ the coefficient matrix of $\fzon$ at time $t_n$, i.e. 
$$
\fzon(t_n,\xvec,\vvec) = e^{-\left|\vvec \right|^2} \sum\limits_{j = 0}^{\ndof_\xvec}\sum\limits_{m = 0}^{\ndof_\vvec}c_{j,m}^{n}\phi_j(\xvec)L_m\left(\vvec\right).
$$
For the parameters $\Vvec_n$ and $T_n$ we let $\Vvec_n := \Vvec(t_n,\,.\,)$ as well as $T_n=T(t_n,\,.\,)$. In order to update $c^{n} \rightarrow c^{n+1}$, we use \refer{eq:variational_formulation_zo_final} and replace $t_0$ by $t_n$ as well as $t_{n_t}$ by $t_{n+1}$. To obtain an explicit scheme, we approximate the remaining $dt$-integrals by a left sided rectangle rule
\begin{equation*}
\int\limits_{t_n}^{t_{n+1}} g(t)\, dt \approx  (t_{n+1} -t_n) g(t_n).
\end{equation*}
Last but not least, we define matrices corresponding to the integrals in \ref{eq:variational_formulation_zo_final} to shorten notation. Dropping the superscripts $\Vvec$ and $T$ from the basis functions, we have
\begin{equation*}
\begin{aligned}
M^n_{i,j} &:= \sum\limits_{\el\in\mesh}\int\limits_\el e^{-|\vvec|^2}T^{\nicefrac{3}{2}} \varphi_j\varphi_i\,d(\xvec,\vvec) & \text{weighted mass matrix} \\
G^n_{ij}&:=\sum\limits_{\el\in\mesh}\int\limits_\el e^{-|\vvec|^2}T^{\nicefrac{3}{2}} \varphi_j\mathbb{D}_t\nabla_\vvec\varphi_i\,d(\xvec,\vvec)&\text{ time derivatives}\\
F^n_{ij}&:=\sum\limits_{\el\in\mesh}\int\limits_{\partial \el}e^{-|\vvec|^2}T^{\nicefrac{3}{2}}(\sqrt{T}\vvec+\Vvec) \cdot \nvec \varphi_j\varphi_i\,d(\xvec,\vvec)&\text{fluxes}\\&-\int\limits_\el e^{-|\vvec|^2}T^{\nicefrac{3}{2}} \varphi_j\mathbb{D}_\xvec\nabla\varphi_i\,d(\xvec,\vvec). &
\end{aligned}
\end{equation*}
The application of the collision operator to a given coefficient vector $c$ is dentod as $Q^n(c)c$. The superscript $n$ is used to highlight the time dependency of the above matrices. This time dependence has to be understood in terms of time dependent $\Vvec$ and $T$. In matrices superscripted with $n$, the ansatz velocity $\Vvec$ and temperature $T$ are evaluated at $t_n$.\par
Considering \refer{eq:variational_formulation_zo_final}, we note that the first term -- which we need to perform the update $c^{n} \rightarrow c^{n+1}$ -- contains $T_{n+1}$, the scheme is not explicit. However, due to the huge numerical work required for the collision integrals, an explicit scheme is likely. We propose the following algorithm to avoid implicit time stepping:
\begin{enumerate}
\item calculate a helper distribution $h$, assuming constant (w.r.t time) $\Vvec_n$ and $T_n$. $h$ is the solution at time $t_{n+1}$, expressed with $\Vvec_n$ and $T_n$.
\item calculate updated values for $(\Vvec_f)_{n+1}$ and $(T_f)_{n+1}$ and smooth them using \refer{eq:variational_formulation_V_and_T} to obtain $\Vvec_{n+1}$ and $T_{n+1}$.
\item advance the distribution in time, taking the evolution of $\Vvec_n$ and $T_n$ into account.
\end{enumerate}
At first sight, it seems we double the numerical work with this strategy by doing the update twice. This can be avoided as we show in the sequel. \\
Let us consider the calculation of $h$:
We perform a step with fixed $\Vvec_n$ and $T_n$ and step size $\tau$. Let us also -- with a slight abuse of notation --  denote the resulting coefficient vector by $h$, which reads
\begin{equation}
h = c^{n} - \tau\left(M^n\right)^{-1} \left( F^n c^{n} - Q^n(c^{n})c^{n}\right).
\label{eq:matrix_formulation_d}
\end{equation}
From the distribution $h$, $\Vvec_n$ and $T_n$ we arrive at $\Vvec_f(t_{n+1},.)$ and $T_f(t_n,.)$ smooth them using \refer{eq:variational_formulation_V_and_T}, and obtain $\Vvec_{n+1}$ and $T_{n+1}$. With these values, we perform the update, taking the time evolution of $\Vvec$ and $T$ into account. The new coefficient vector for time $t^{n+1}$ satisfies
\begin{equation*}
\begin{aligned}
M^{n+1} c^{n+1} - M^{n} c^{n} + \tau(F^n-G^n)c^{n} -\tau Q^n(c^{n})c^{n} = 0.
\end{aligned}
\end{equation*}
Solving for the unknown $c^{(n+1)}$, we have
\begin{equation}
\begin{aligned}
c^{n+1} = \left(M^{n+1}\right)^{-1} \left(M^{n} c^{n} + \tau(G^n-F^n)c^{n} + \tau Q^n(c^{n}) c^{n} \right).
\end{aligned}
\label{eq:matrix_formulation}
\end{equation}
Finally, we note that \refer{eq:matrix_formulation_d} and \refer{eq:matrix_formulation} can be combined to simplify the calculation of $c^{n+1}$
$$
c^{n+1} = \left(M^{n+1}\right)^{-1} \left(M^{n} h + \tau G^n c^{n}\right).
$$
\subsection{\texorpdfstring{Choosing $\Vvec$ and $T$}{ChoosingVandT}}
The choice of the ansatz velocity and the ansatz temperature are very crucial to improve stability. We expect optimal approximation properties if $\Vvec = \Vvec_f$ and $T = T_f$. However, that results in constraints to the approximation space, from proposition \ref{prop:macroscopics_of_fzo} we deduce that $\Vzo \equiv 0$ and $\Tzo \equiv 1$ have to hold. Additionally, when we consider shocks, high order solutions are subject to spurious osscillations. These are fed back into the transformed equation. In order to get rid of these constraints and also to improve stability we use smoothed values for $\Vvec$ and $T$. Let us also briefly discuss the boundary conditions: On the inflow boundary we use stationary dirichlet values for both, $\Vvec$ and $T$. Along obstacles we use consistent boundary conditions for $\Vvec$: $\Vvec \cdot \nvec \equiv \mathbf{0}, \xvec \in (\partial\Omega)_\text{spec}$ and $\Vvec \equiv \mathbf{0}$ respectively. For $T$ we use natural boundary conditions along obstacles. Thus, we consider the following $H_0^1$ problems for $\Vvec$ and $T$:
\begin{equation}
\begin{aligned}
\text{Find } &T \in H_0^1(\Omega),\;\, T =T_f(t = 0),\,\xvec\in(\partial\Omega)_\text{in} \text{ s.t.:}\\
&\int\limits_\Omega T v\,d\xvec + \lambda \int\limits_\Omega \nabla T \nabla v\,d\xvec = \int\limits_\Omega {T_f} v\,d\xvec & \forall\;v \in H_0^1(\Omega) \\
\text{Find } &\Vvec \in [H_0^1(\Omega)]^2,\;\, \Vvec = \Vvec_f(t = 0),\,\xvec\in(\partial\Omega)_\text{in}\cup(\partial\Omega)_\text{diff} \text{ s.t.:} \hfill\\
&\int\limits_\Omega \Vvec v\,d\xvec + \lambda \int\limits_\Omega \nabla\Vvec \nabla v\,d\xvec+\frac{1}{\epsilon}\!\!\!\int\limits_{(\partial\Omega)_{\text{spec}}}\!\!\!(\Vvec\cdot \nvec) (v\cdot \nvec)\,d\xvec = \int\limits_\Omega ({\Vvec_f}) v \,d\xvec& \forall\;v \in [H_0^1(\Omega)]^2\\
\end{aligned}
\label{eq:variational_formulation_V_and_T}
\end{equation}
The parameter $\lambda$ is chosen according to the polynomial degree w.r.t. $\xvec$ and the mesh size $h$, i.e. $\lambda = c\tfrac{h^2}{p^2}$. Specific values for the parameter $c$ are stated in section \ref{sec:numerical_results}. 
\begin{remark}
The rescaling idea was already used in two velocity dimensions in the PHD thesis of the first author, \cite{Kitzler2016}. Therein, we used piecewise constant (w.r.t. time and physical space) quantities $\Vvec$ and $T$, instead of continuous ones. Consequently, a time step was composed of two sub steps: First we advance the density and calculate updated values for $\Vvec$ and $T$. In a second step, we project the coefficient vector (which represents the updated $f$ with the old values of $\Vvec$ and $T$) to the updated values of $\Vvec$ and $T$. \newline In \cite{FILBET2013177}, the authors apply the rescaling idea to the fast Fourier spectral method as presented in \cite{filbetnlogn}, and apply it to a kinetic equation in granular media. They also treat the specific case of the Boltzmann equation, where the scaling function was made independent of $f$, using the conservation laws \refer{eq:unclosed_cons_law} and a Maxwell closure. In contrast to our pull back which is based on temperature and velocity, they use a scaling function based on the internal energy. \\
In \cite{20.500.11850/152993,Kllermeier:698099}, rescaling is used for the numerical solution of the BGK equation in a finite volume method. The velocity discretization is as it is here, a Maxwellian weighted, polynomial spectral method. 
\end{remark}
\subsection{Application of the collision operator}
The application of the collision operator \refer{eq:collop} results typically in computational complexity $\mathcal{O}(\ndof^3)$. We developed an algorithm which has a lower complexity $\mathcal{O}(\ndof^{\nicefrac{7}{3}})$ to reduce computational costs, \cite{Kitzler2018}. To illustrate the main ideas, we consider a non negative density distribution, $f = f(\vvec)$. The algorithm is based on 
\begin{itemize}
\item \textbf{Transformation to mean and relative velocities.} We let $\vbar:= \frac{\vvec+\wvec}{2},\,\hat{\vvec}:= \frac{\vvec-\wvec}{2}$ and use a well known representation of the collision operator \cite{Cercignani_1994} to write
\begin{equation*}
\begin{aligned}
&\intrd Q(f)(\vvec)\phi(\vvec)\,d\vvec \\&= 8\intrd \intrd \int\limits_{S^2}b_r(2|\hat{\vvec}|)b_\theta(\tfrac{\hat{\vvec}\cdot \evec'}{|\hat{\vvec}|}) f(\vbar+\hat{\vvec})f(\vbar-\hat{\vvec})[\phi(\vbar+\evec'|\hat{\vvec}|)-\phi(\vbar+\hat{\vvec})]\,d\evec'd\hat{\vvec}\,d\vbar.
\end{aligned}
\end{equation*}
We interpolate the product of $f$ evaluations in $V_{2N}$, what is efficient due to the Lagrange polynomials unity Vandermonde matrix.
\item \textbf{basis polynomials, diagonalizing the inner integrals.} We provide a polynomial basis, consisting of Spherical harmonics on the unit sphere and of generalized Laguerre polynomials in the radial direction. This basis diagonalizes the two innermost integrals w.r.t. to $\vhat$ and $\evec'$. The $\vbar$ integral is computed by a Gauss Hermite quadrature rule.
\item \textbf{efficient basis transformations.} To transform the solution from nodal representation to the spherical one, the transform is
\begin{equation*}
\mathbf{L} \rightarrow \mathbf{H} \rightarrow \Theta\rightarrow \mathbf{\Psi},
\end{equation*}
where $\mathbf{L}$ denotes the Lagrange representation of $\fzo$. $\mathbf{H}$ is $\fzo$ in the Hermite polynomial basis, an efficient transformation is obtained, using the tensor product structure in both of the bases, factorizing the three dimensional transformation into several one dimensional ones.
The second and third part of the transformation are efficient, since all of the involved polynomial bases provide hierarchical orthogonal decompositions resulting in sparse transformations.
\par
For more details, we refer the reader to \cite{Kitzler2018}.
\end{itemize}
\section{Numerical results}
The results presented in this section were obtained using the finite element library NgSolve \cite{ngsolve}.
\label{sec:numerical_results}
\subsection{Shock tube}
\label{sec:numerical_results1}
\begin{figure}
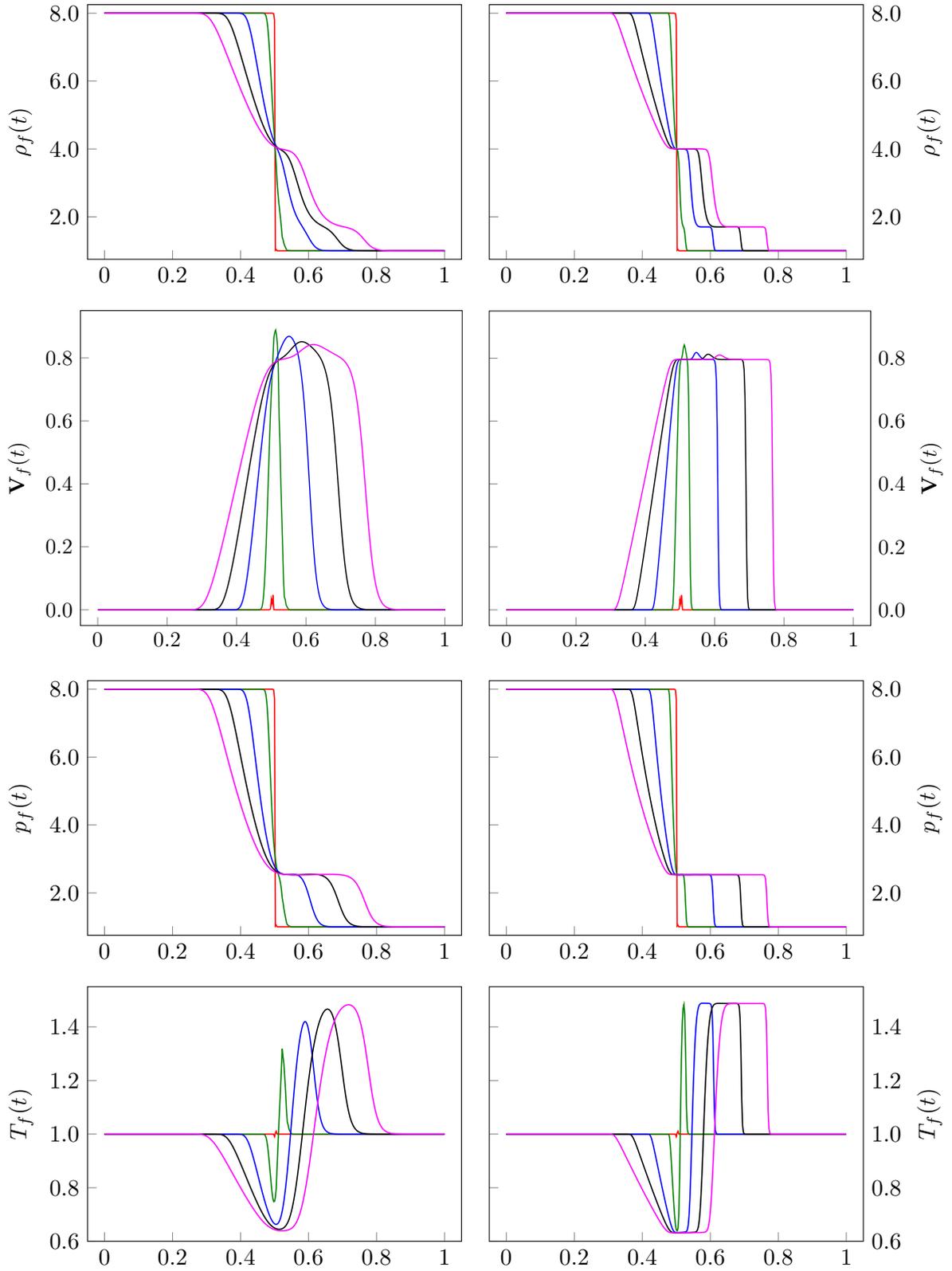

    \centering
    \begin{subfigure}[b]{0.49\textwidth}
		\resizebox{\linewidth}{!}{\input{./images/shock_tube/time_evolution/rho_200_els_kn0.01_artvs_2.5e-05.tikz}}
    \end{subfigure}
    \hfill
    \begin{subfigure}[b]{0.49\textwidth}
		\resizebox{\linewidth}{!}{\input{./images/shock_tube/time_evolution/rho_200_els_kn0.001_artvs_2.5e-05.tikz}}
    \end{subfigure}

    \vspace{0.3cm}

    \begin{subfigure}[b]{0.49\textwidth}
		\resizebox{\linewidth}{!}{\input{./images/shock_tube/time_evolution/v_200_els_kn0.01_artvs_2.5e-05.tikz}}
    \end{subfigure}
    \hfill
    \begin{subfigure}[b]{0.49\textwidth}
		\resizebox{0.999\linewidth}{!}{\input{./images/shock_tube/time_evolution/v_200_els_kn0.001_artvs_2.5e-05.tikz}}
    \end{subfigure}

    \vspace{0.3cm}
        
    \begin{subfigure}[b]{0.49\textwidth}
		\resizebox{\linewidth}{!}{\input{./images/shock_tube/time_evolution/pressure_200_els_kn0.01_artvs_2.5e-05.tikz}}
    \end{subfigure}
    \hfill
    \begin{subfigure}[b]{0.49\textwidth}
		\resizebox{\linewidth}{!}{\input{./images/shock_tube/time_evolution/pressure_200_els_kn0.001_artvs_2.5e-05.tikz}}
    \end{subfigure}

    \vspace{0.3cm}
        
    \begin{subfigure}[b]{0.49\textwidth}
		\resizebox{\linewidth}{!}{\input{./images/shock_tube/time_evolution/temperature_200_els_kn0.01_artvs_2.5e-05.tikz}}
    \end{subfigure}
    \hfill
    \begin{subfigure}[b]{0.49\textwidth}
		\resizebox{\linewidth}{!}{\input{./images/shock_tube/time_evolution/temperature_200_els_kn0.001_artvs_2.5e-05.tikz}}
    \end{subfigure}        
    
	\caption{Snapshots at $t = dt, 0.014, 0.056,0.098,0.14$ for the shock tube example in section \ref{sec:numerical_results1}. The calculation was done with order 9 polynomials in $\vvec$ and order 4 polynomials in $ \xvec$. The knudsen numbers are $\text{kn}=0.01$ (left) and $\text{kn}=0.001$ (right) respectively.}
	\label{fig:shocktube_timeevo}

\end{figure}

The first test we perform is a shock tube on the domain $\R{} \times \R3$. We split the $\xvec$ axis into the positive and the negative half, separated by a diaphragm which is suddenly removed at time $t = 0$. The two initial states of the gas are assumed to be locally equilibrated:
$$
f(0,\xvec,\vvec) = \frac{\rho_f(\xvec)}{(2\pi T_f(\xvec))^{\nicefrac{3}{2}}}e^{-\left|\frac{\vvec-\Vvec_f}{\sqrt{2T_f}} \right|^2}.
$$
The values for the initial density, velocity and temperature are
\begin{equation}
\rho_f(\xvec) = \begin{cases}8 & \xvec \leq 0 \\ 1 & \text{else}\end{cases}, \quad\quad
T_f(\xvec) \equiv 1, \quad\quad
\Vvec_f(\xvec) \equiv \mathbf{0}.
\end{equation}
The time evolution for the above problem is shown in figure \ref{fig:shocktube_timeevo} for knudsen numbers $\text{kn} = 0.01$ and $\text{kn} = 0.001$ respectively. The results were obtained with order 9 polynomials for $\vvec$ and order 4 polynomials in $\xvec$, the mesh for the calculations consists of 200 aequidistant elements only. Time stepping was done with the explicit runge kutta 4 scheme, using a very small step $\tau = 1.5625 \text{e-}5$, such that we can neglect the temporal error.
In figure \ref{fig:shocktube_kn001} we compare different polynomial orders w.r.t. $\vvec$ and also different values of $c$ in \refer{eq:variational_formulation_V_and_T}. In the first row there's hardly any difference to see, pressure and density are independent of the scaling parameter. Turning to the second row, we note that bulk velocity and temperature are essentially independent of the constant $c$ also. Considering the polynomial degree $N$, we note a small difference between the solutions for $N=3$ and $N = 4$. The solutions for $N=4$ and $N=5$ coincide perfectly. In the bottom row the shift and scale parameters $\Vvec$ and $T$ are shown. Clearly, both of them depend significantly on the viscosity parameter. In contrast, the dependence on the polynomial degree $N$ is rather weak.

\subsection{Flow around airfoil}
\label{sec:numerical_results2}
\begin{figure}
    \centering
		\resizebox{0.47\linewidth}{!}{\includegraphics[trim=0.0px 0.0px 800px 0px,clip]{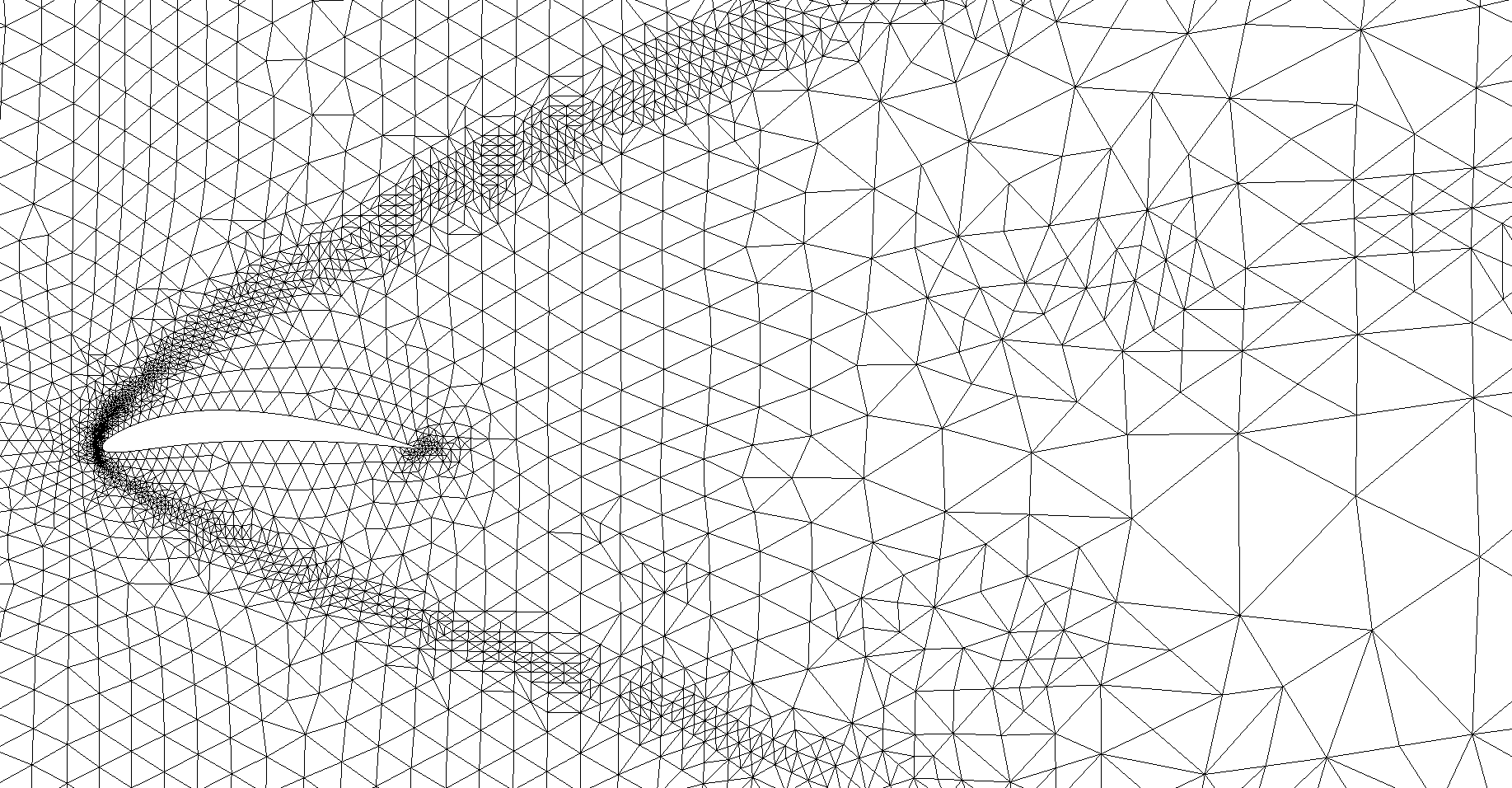}}
	\hfill
		\resizebox{0.47\linewidth}{!}{\includegraphics[trim=0.0px 0.0px 800px 0px,clip]{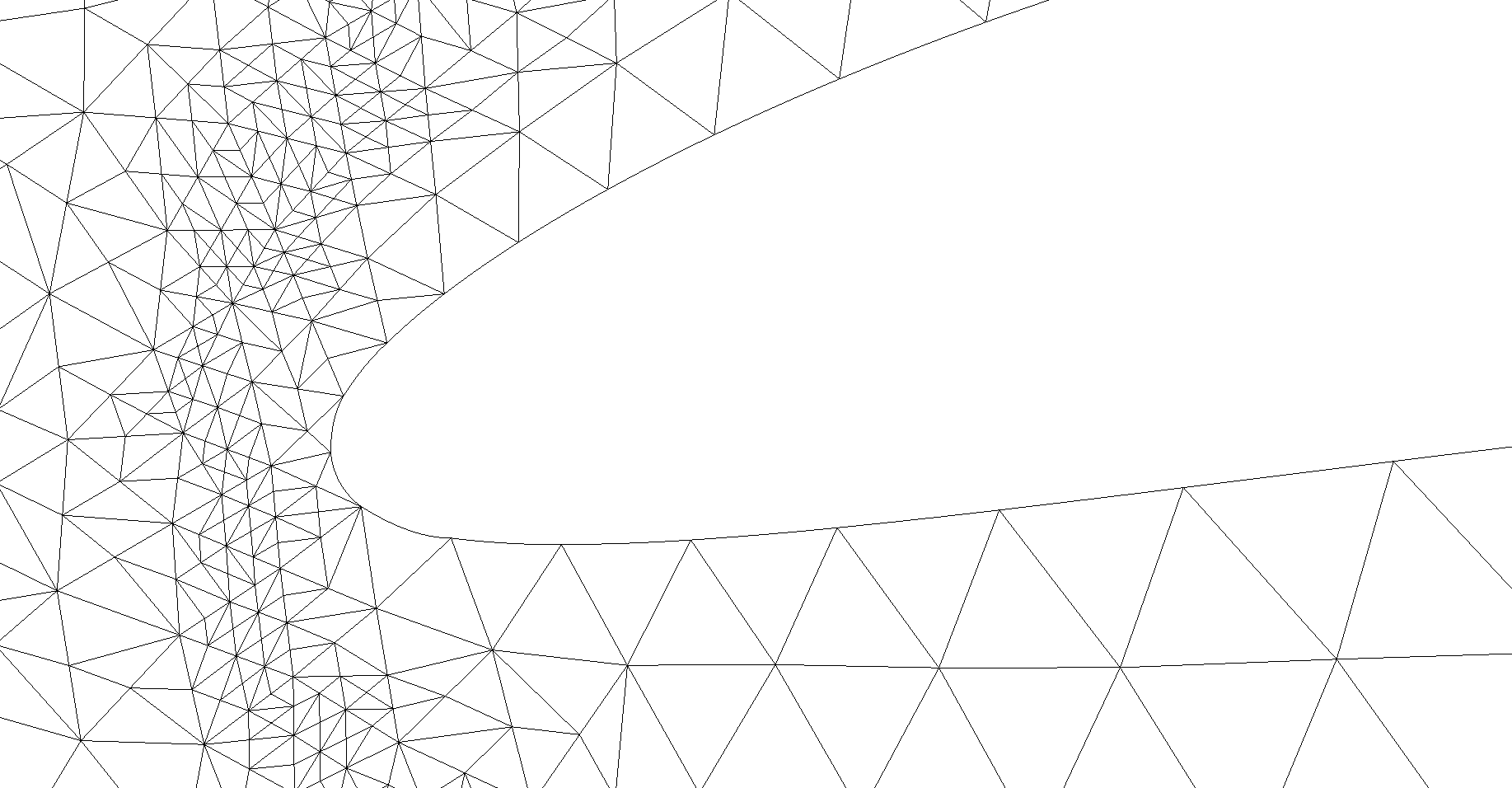}}
	\caption{Computational mesh for the example in section \ref{sec:numerical_results2}.}
	\label{fig:NACA7410_mesh}
\end{figure}

\begin{figure}[htb]
\centering
\begin{tikzpicture}
\begin{axis}[
hide axis,
scale only axis,
height = 0pt,
width = 0pt,
colormap name = ngsolve,
colorbar sampled,
colorbar horizontal,
point meta min = 0.0,
point meta max = 875,
colorbar style ={samples = 9,width = 10cm, height = 0.3cm, xtick={0.0,437.5,875}, xticklabel style = {style={font=\footnotesize}},
xticklabel pos = top,
        xticklabel={
    \pgfkeys{/pgf/fpu=true, /pgf/fpu/output format=fixed}
    \pgfmathparse{\tick*1e-5}
            \num[
                scientific-notation = false,
                exponent-product=\cdot,
                output-exponent-marker = \text{e},
            ]{\pgfmathresult}
        }
}]
\addplot [draw = none] coordinates { (0,0) };
\end{axis}
\end{tikzpicture} 

\vspace{0.25cm}
    \begin{subfigure}[b]{0.49\textwidth}
		\resizebox{0.99\linewidth}{!}{\includegraphics[trim=0.0px 130.0px 0px 170px,clip]{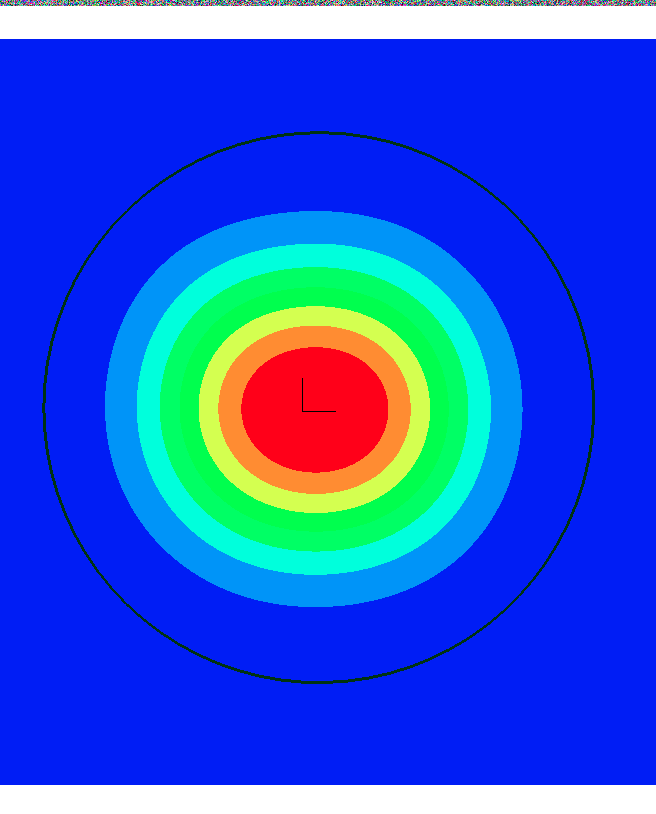}}
	\subcaption{$(x,y) = (-0.5,0.008)$}
	\label{fig:distribution_bow_1}	
    \end{subfigure}
    \hfill    
    \begin{subfigure}[b]{0.49\textwidth}
		\resizebox{0.99\linewidth}{!}{\includegraphics[trim=0.0px 130.0px 0px 170px,clip]{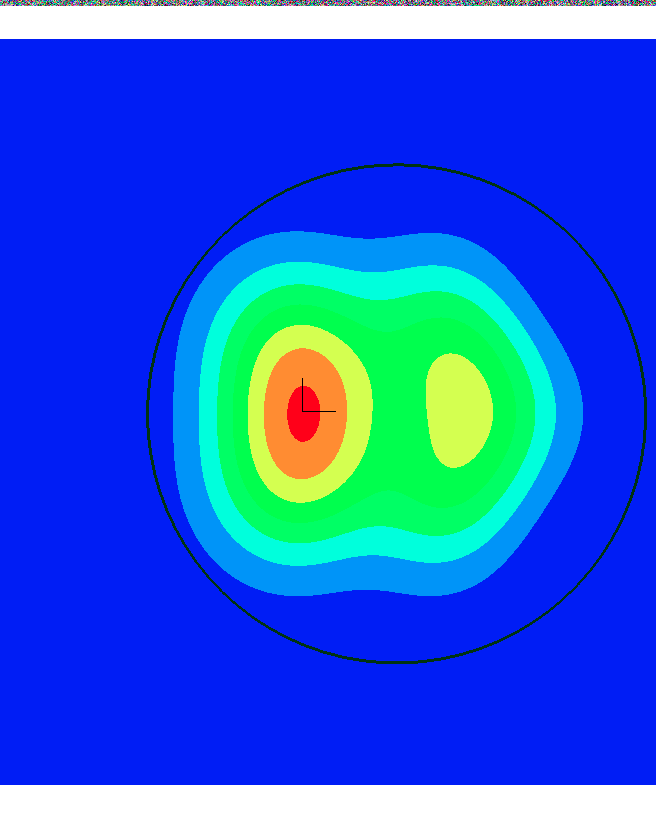}}
	\subcaption{$(x,y) = (-0.511875,0.008)$}
	\label{fig:distribution_bow_2}
    \end{subfigure}
    
    \vspace{0.25cm}
    
    \centering
    \begin{subfigure}[b]{0.49\textwidth}
		\resizebox{0.99\linewidth}{!}{\includegraphics[trim=0.0px 130.0px 0px 170px,clip]{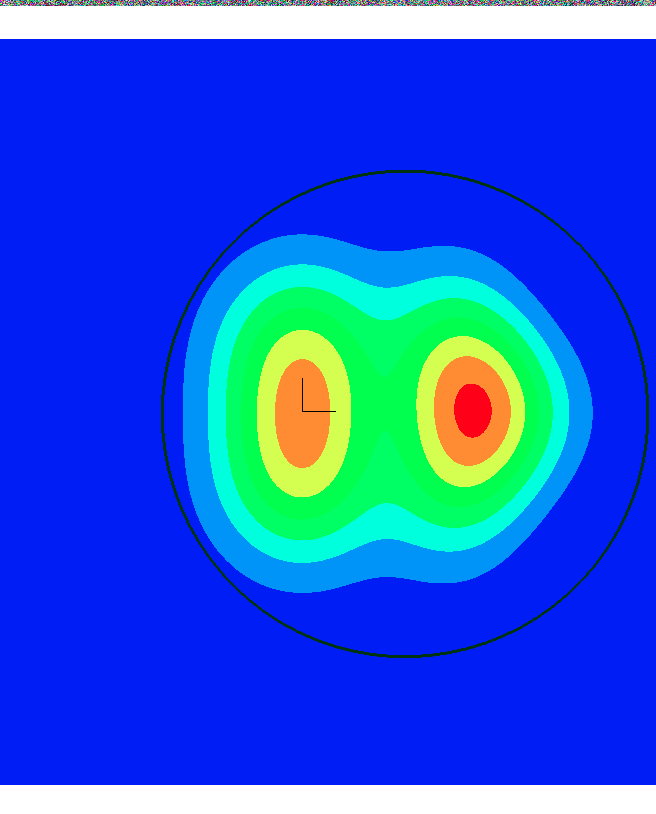}}
	\subcaption{$(x,y) = (-0513125,0.008)$}
	\label{fig:distribution_bow_3}	
    \end{subfigure}
    \hfill    
    \begin{subfigure}[b]{0.49\textwidth}
		\resizebox{0.99\linewidth}{!}{\includegraphics[trim=0.0px 130.0px 0px 170px,clip]{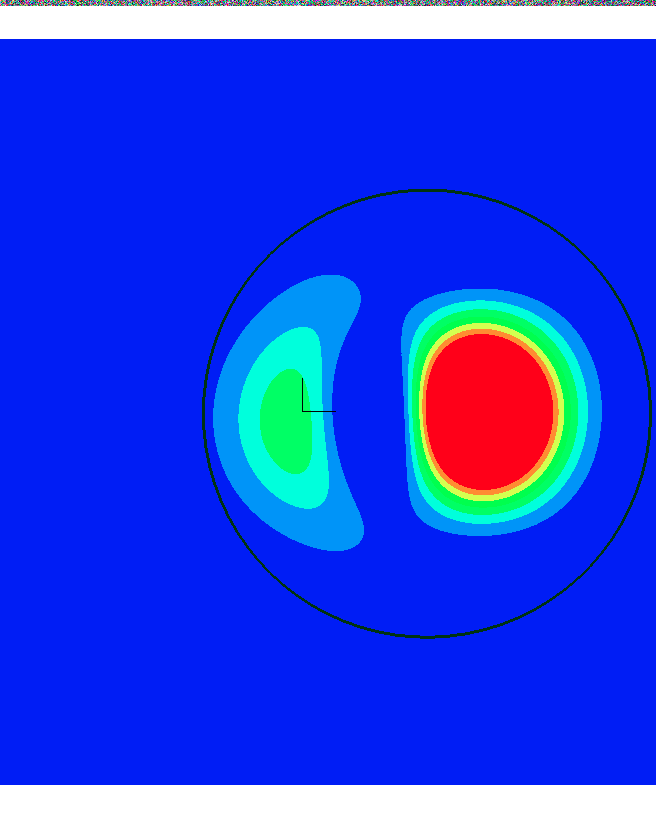}}
	\subcaption{$(x,y) = (-0.51625,0.008)$}
	\label{fig:distribution_bow_4}
    \end{subfigure}    
    
    
    \centering
    \begin{subfigure}[b]{0.49\textwidth}
		\resizebox{0.99\linewidth}{!}{\includegraphics[trim=0.0px 130.0px 0px 170px,clip]{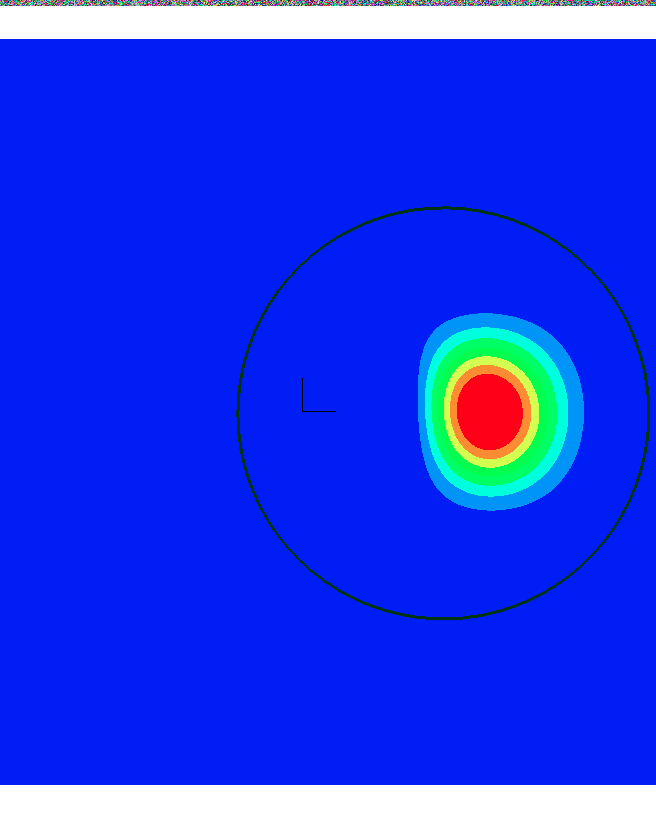}}
	\subcaption{$(x,y) = (-0.51875,0.008)$}
	\label{fig:distribution_bow_5}	
    \end{subfigure}
    \hspace{0.5cm}
\begin{tikzpicture}
\begin{axis}[
hide axis,
scale only axis,
height = 0pt,
width = 0pt,
colormap name = ngsolve,
colorbar sampled,
point meta min = 0.0,
point meta max = 25,
colorbar style ={samples = 9,width = 0.3cm, height = 6.4cm, ytick={0.0,12.5,25}, xticklabel style = {style={font=\footnotesize}},
yticklabel pos = right,
        yticklabel={
    \pgfkeys{/pgf/fpu=true, /pgf/fpu/output format=fixed}
    \pgfmathparse{\tick*1e-3}
            \num[
                scientific-notation = false,
                exponent-product=\cdot,
                output-exponent-marker = \text{e},
            ]{\pgfmathresult}
        }
}]
\addplot [draw = none] coordinates { (0,0) };
\end{axis}
\end{tikzpicture}
%
    \caption{The density at several positions through the bow shock. 
The upper color bar is for \subref{fig:distribution_bow_1}-\subref{fig:distribution_bow_4}, the vertical bar is for \subref{fig:distribution_bow_5}.}
\label{fig:distributions_bow_shock}
\end{figure}

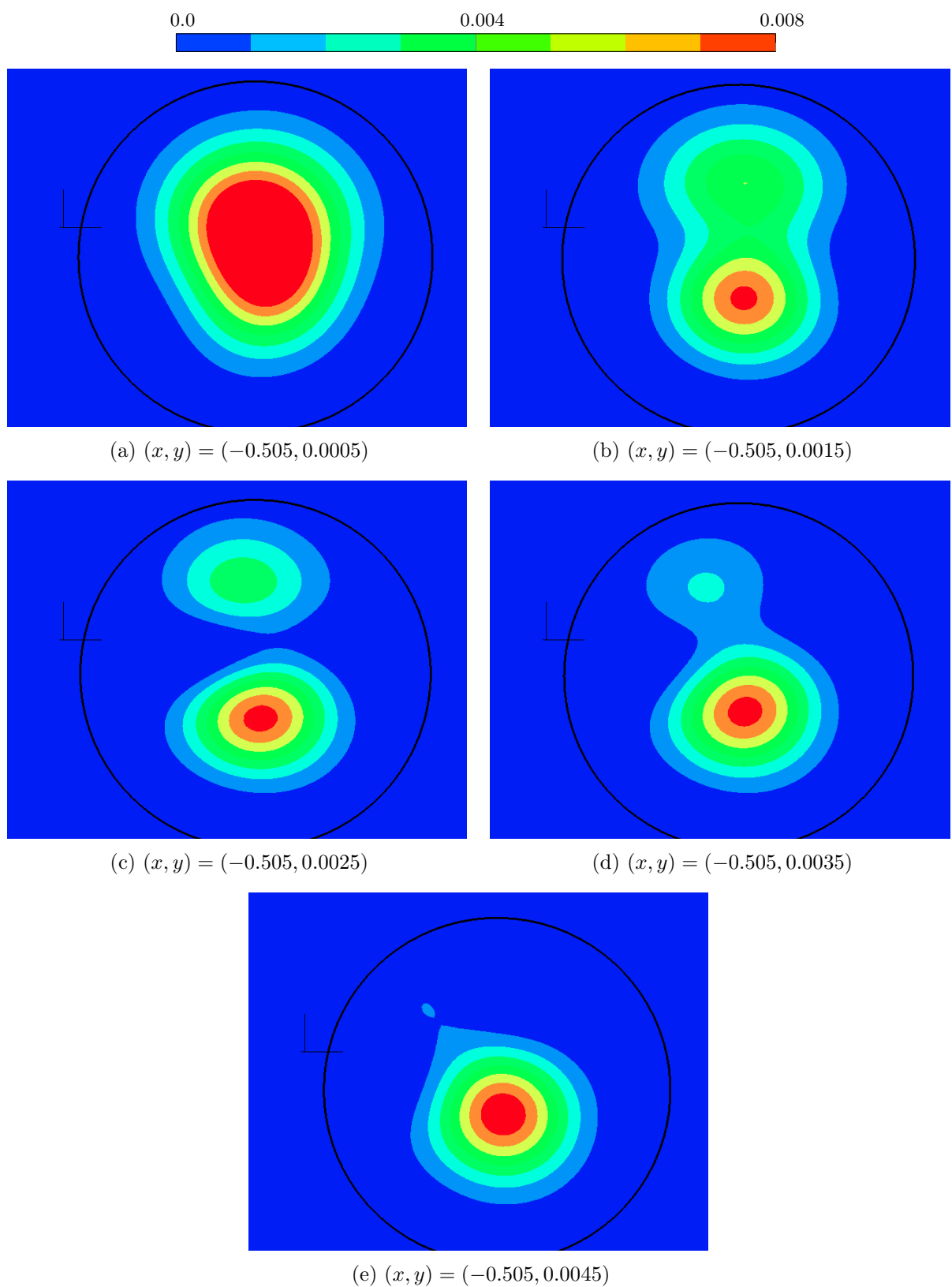
\begin{figure}[htb]
    \centering
\begin{tikzpicture}
\begin{axis}[
hide axis,
scale only axis,
height = 0pt,
width = 0pt,
colormap name = ngsolve,
colorbar sampled,
colorbar horizontal,
point meta min = 0.0,
point meta max = 8,
colorbar style ={samples = 9,width = 10cm, height = 0.3cm, xtick={0.0,4,8}, xticklabel style = {style={font=\footnotesize}},
xticklabel pos = top,
        xticklabel={
    \pgfkeys{/pgf/fpu=true, /pgf/fpu/output format=fixed}
    \pgfmathparse{\tick*1e-3}
            \num[
                scientific-notation = false,
                exponent-product=\cdot,
                output-exponent-marker = \text{e},
            ]{\pgfmathresult}
        }
}]
\addplot [draw = none] coordinates { (0,0) };
\end{axis}
\end{tikzpicture}
   
   \vspace{0.25cm}
    \begin{subfigure}[b]{0.49\textwidth}
		\resizebox{0.99\linewidth}{!}{\includegraphics[trim=0.0px 130.0px 0px 170px,clip]{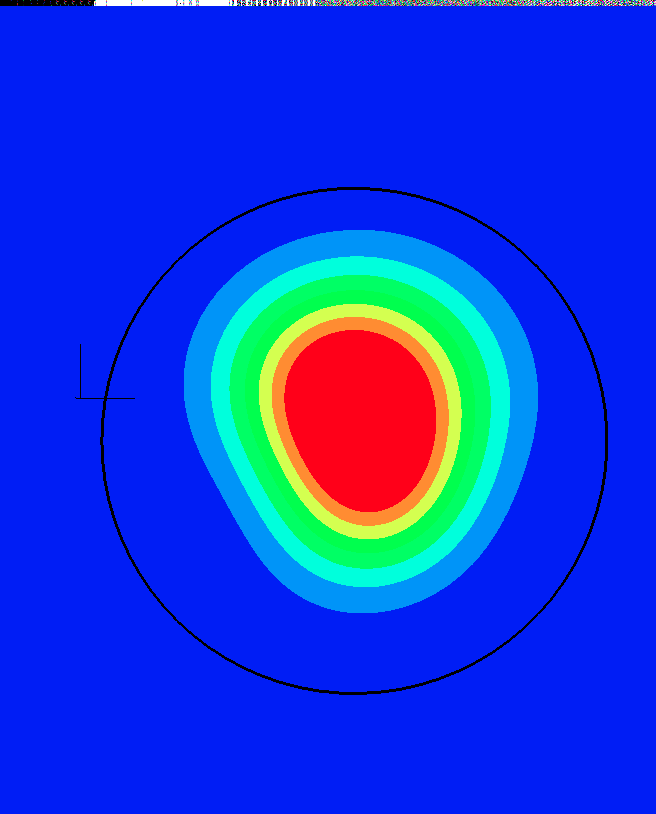}}
	\subcaption{$(x,y) = (-0.505,0.0005)$}
	\label{fig:distribution_rear_1}	
    \end{subfigure}
    \hfill    
    \begin{subfigure}[b]{0.49\textwidth}
		\resizebox{0.99\linewidth}{!}{\includegraphics[trim=0.0px 130.0px 0px 170px,clip]{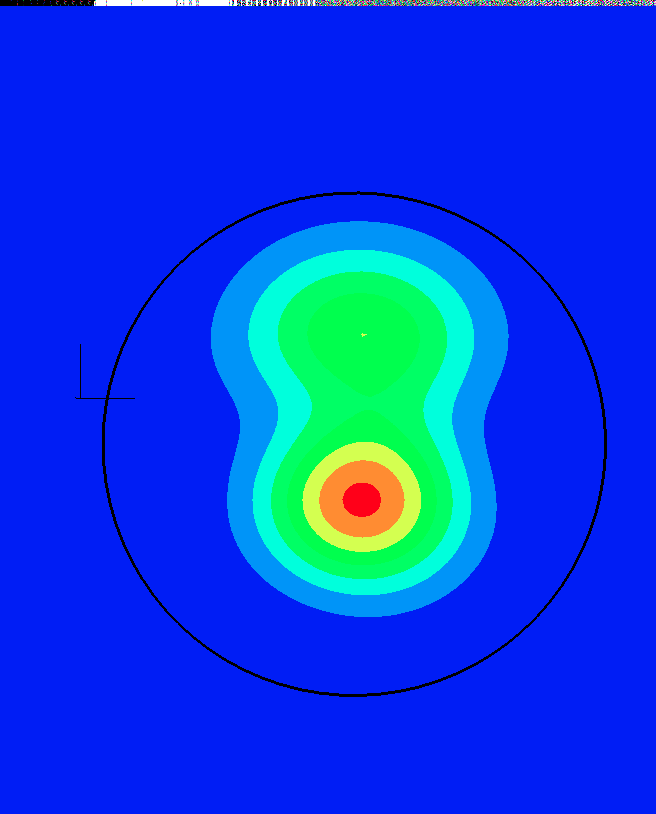}}
	\subcaption{$(x,y) = (-0.505,0.0015)$}
	\label{fig:distribution_rear_2}
    \end{subfigure}

    \vspace{0.25cm}  
    \centering
    \begin{subfigure}[b]{0.49\textwidth}
		\resizebox{0.99\linewidth}{!}{\includegraphics[trim=0.0px 130.0px 0px 170px,clip]{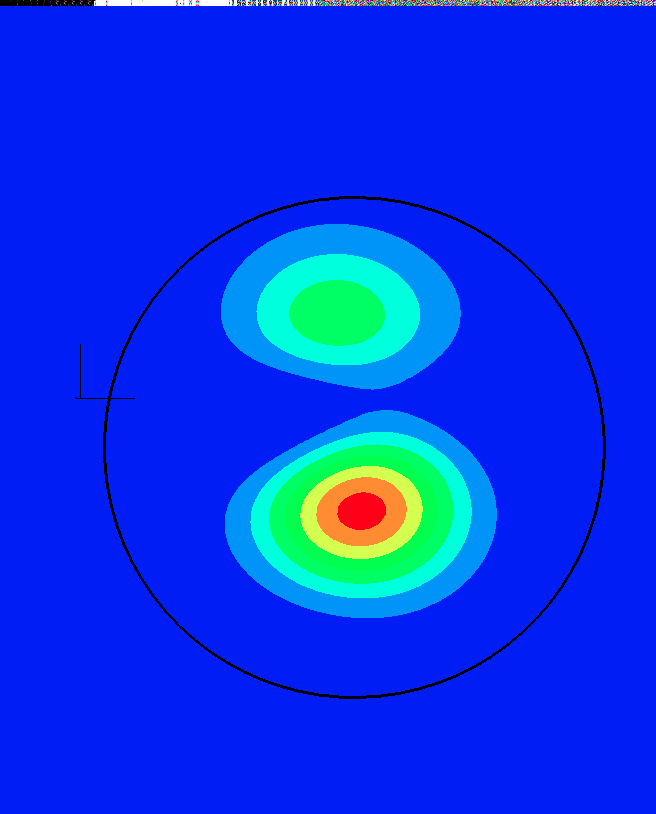}}
	\subcaption{$(x,y) = (-0.505,0.0025)$}
	\label{fig:distribution_rear_3}	
    \end{subfigure}
    \hfill    
    \begin{subfigure}[b]{0.49\textwidth}
		\resizebox{0.99\linewidth}{!}{\includegraphics[trim=0.0px 130.0px 0px 170px,clip]{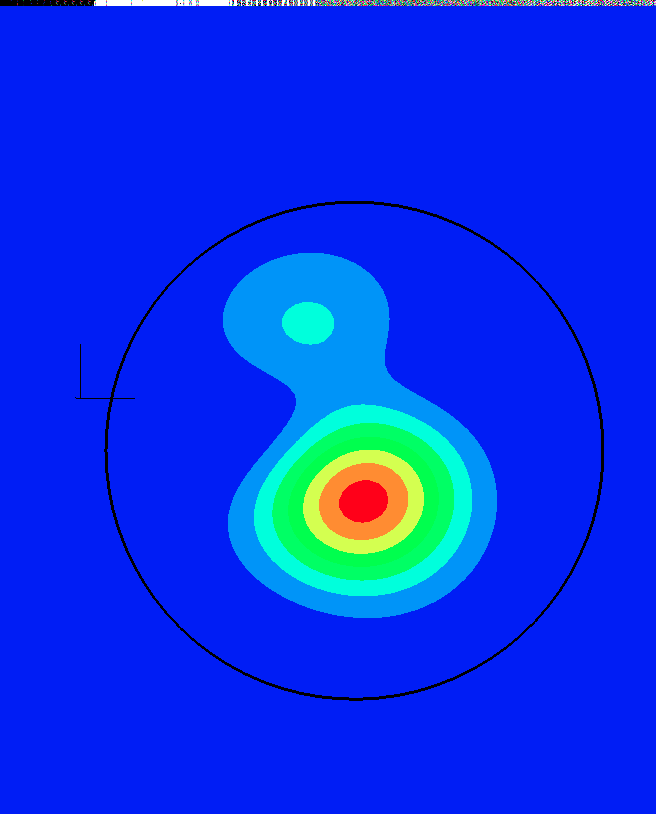}}
	\subcaption{$(x,y) = (-0.505,0.0035)$}
	\label{fig:distribution_rear_4}
    \end{subfigure}    

    \vspace{0.25cm}  
    \begin{subfigure}[b]{0.49\textwidth}
		\resizebox{0.99\linewidth}{!}{\includegraphics[trim=0.0px 130.0px 0px 170px,clip]{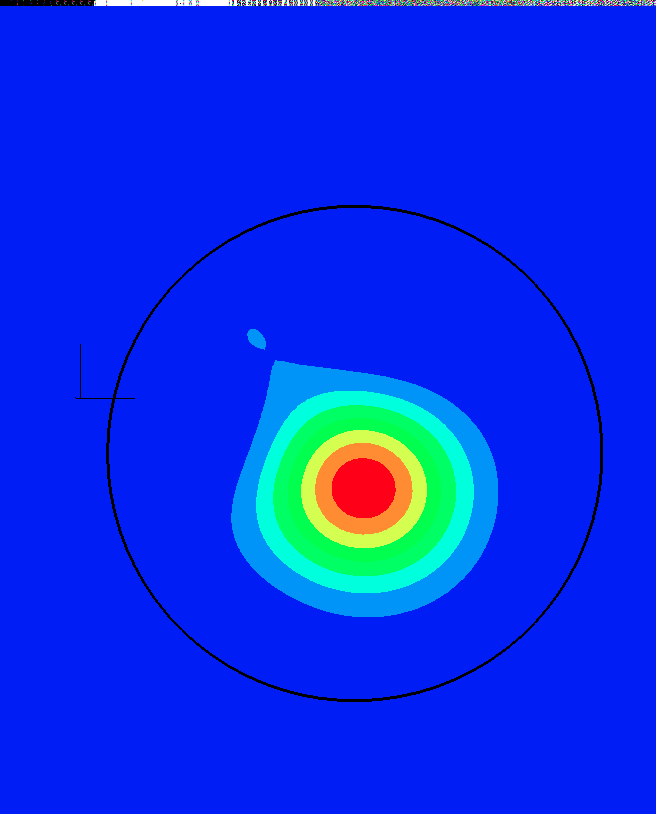}}
	\subcaption{$(x,y) = (-0.505,0.0045)$}
	\label{fig:distribution_rear_5}	
    \end{subfigure}
    \caption{The density through the trailing edge shock.}
    \label{fig:distributions_rear_shock}
\end{figure}

In the second example we consider the flow around a NACA 7410 air foil, the free flow mach number being $m_{\infty} = 4.5$. As in the previous example, we assume that the gas is initially in local equilibrium. To obtain the desired mach number, we let
\begin{equation*}
f(0,\xvec,\vvec) = \frac{\rho_f(\xvec)}{(2\pi T_f(\xvec))^{\nicefrac{3}{2}}}e^{-\left|\frac{\vvec-\Vvec_f}{\sqrt{2T_f}} \right|^2}, \quad\quad \text{ where } \rho_f(\xvec) \equiv 1, \quad T_f(\xvec) \equiv 1.
\end{equation*}
The initial velocity profile is constant almost everywhere (with approximate value 5.81), and additionally satisfies $\Vvec\cdot \nvec \equiv 0$ along the air foil. This velocity profile was obtained by solving \refer{eq:variational_formulation_V_and_T}, with $\lambda_0 = 2.5\text{e-}3$, the penalty term for vanishing normal velocity is rather small, $\tfrac{1}{\epsilon}=100$. During time stepping -- which was a forward euler scheme using $\tau = 7.8125\text{e-}6$ -- the penalty parameter was kept at its initial value, $\lambda$ was chosen as in \refer{eq:variational_formulation_V_and_T} with $c = 18$.\\
For the results presented in figures \ref{fig:NACA7410} - \ref{fig:distributions_rear_shock} we chose a Knudsen number $\text{kn} = 0.005$. The polynomial degrees are 6 for $\xvec$ and 3 for $\vvec$. The mesh for this computation consisted of 6862 elements, part of it is shown in figure \ref{fig:NACA7410_mesh}. It was obtained by adaptive refinement along the shocks, using the ZZ-error estimator developed in \cite{NME:NME1620330702,NME:NME1620330703}.
Figure \ref{fig:NACA7410} shows the macroscopic behaviour in a large domain around the air foil, the right column shows details of the bow shock. Note that there are no capturing techniques applied, yet hardly any oscillations nor overshoots are observed, we obtain a smooth shock front.   \linebreak Figures \ref{fig:distributions_bow_shock} and \ref{fig:distributions_rear_shock} show distributions throguh the bow shock and along the trailing edge, which are far from local equilibrium. Even though only 64 basis functions in $\vvec$ are used in this calculation, the distributions are well resolved. Figures \ref{fig:distributions_bow_shock} and \ref{fig:distributions_rear_shock} also justify the usage of a Boltzmann solver, hydrodynamic models are not suitable at the spatial points depicted in figures \ref{fig:distributions_bow_shock} and \ref{fig:distributions_rear_shock}. However, they also indicate the need for a coupling of macroscopic and mesoscopic equations, in a large part of the domain the distribution is in local equilibrium. \linebreak
This example emphasizes the capabilities of the velocity transformation in combination with a problem suited discretization. Even with a low order expansion, large variations in the mach number can be resolved.

\section{Conclusion}
In this paper we presented a numerical method for the Boltzmann equation. Our focus -- on the one hand -- was to obtain good approximation properties uniform in the mach number. To that end, we proposed to shift and scale the velocity variable by the (smoothed) bulk velocity and temperature. This gave rise to a new distribution $\fzo$ for which we derived the discrete evolution equation. The distribution $\fzo$ has almost everywhere a vanishing bulk velocity and a temperature close to one. That allowed us to use reference trial functions.

To approximate $\fzo$ we proposed Maxwellian weighted polynomials as the trial functions for the velocity. The results in section \ref{sec:numerical_results} confirmed their excellent approximation properties, less than one hundred degrees of freedom in $\vvec$ gave good results, especially for functions close to equilibrium. Polynomial test functions gave the correct conservation laws. 

Finally, the choice of the scaling paramters is a delicate question, as we have seen in numerous numerical experiments. We proposed to extract $\Vvec$ and $T$ from the distribution,and to smooth them with the reaction - diffusion equation \refer{eq:variational_formulation_V_and_T}. Smoothing was required to eliminate spurious osscilations from these quantities. Our results indicate that the solution is almost independent of the specific choice of the smoothing parameter.

\section{Acknowledgments}
Gerhard Kitzler is funded by the Austrian Science Fund (FWF) project F 65.

\begin{figure}
    \centering
    \begin{subfigure}[b]{0.49\textwidth}
		\resizebox{0.999\linewidth}{!}{\input{./images/shock_tube/compare_art_vs/rho_100_els_kn0.01_artvs_2.5e-05.tikz}}
    \end{subfigure}
    \hfill
    \begin{subfigure}[b]{0.49\textwidth}
		\resizebox{0.999\linewidth}{!}{\input{./images/shock_tube/compare_art_vs/pressure_100_els_kn0.01_artvs_2.5e-05.tikz}}
    \end{subfigure}

	\vspace{0.5cm}    
    
    \begin{subfigure}[b]{0.49\textwidth}
		\resizebox{0.999\linewidth}{!}{\input{./images/shock_tube/compare_art_vs/temperature_100_els_kn0.01_artvs_2.5e-05.tikz}}
    \end{subfigure}
    \hfill
    \begin{subfigure}[b]{0.49\textwidth}
		\resizebox{0.999\linewidth}{!}{\input{./images/shock_tube/compare_art_vs/v_100_els_kn0.01_artvs_2.5e-05.tikz}}
    \end{subfigure}

	\vspace{0.5cm}    
        
    \begin{subfigure}[b]{0.49\textwidth}
		\resizebox{0.999\linewidth}{!}{\input{./images/shock_tube/compare_art_vs/temperatureansatz_100_els_kn0.01_artvs_2.5e-05.tikz}}
    \end{subfigure}
    \hfill
    \begin{subfigure}[b]{0.49\textwidth}
		\resizebox{0.999\linewidth}{!}{\input{./images/shock_tube/compare_art_vs/vansatz_100_els_kn0.01_artvs_2.5e-05.tikz}}
    \end{subfigure}
   	\begin{tikzpicture}

\definecolor{color2}{rgb}{0.172549019607843,0.627450980392157,0.172549019607843}
\definecolor{color1}{rgb}{1,0.498039215686275,0.0549019607843137}
\definecolor{color4}{rgb}{0.580392156862745,0.403921568627451,0.741176470588235}
\definecolor{color0}{rgb}{0.12156862745098,0.466666666666667,0.705882352941177}
\definecolor{color3}{rgb}{0.83921568627451,0.152941176470588,0.156862745098039}

\begin{axis}
[
hide axis,
xmin=-2.995, xmax=62.895,
ymin=5.53988011611066e-19, ymax=0.216926593216958,
ymode=log,
tick align=outside,
tick pos=left,
x grid style={lightgray!92.026143790849673!black},
legend cell align={left},
y grid style={lightgray!92.026143790849673!black},
legend columns=3,
legend entries={{$c=64,N=3\quad\quad$},{$c=16,N=3\quad\quad$},{$c=4,N=3\quad\quad$},{$c=64,N=4\quad\quad$},{$c=16,N=4\quad\quad$},{$c=4,N=4\quad\quad$},{$c=64,N=5\quad\quad$},{$c=16,N=5\quad\quad$},{$c=4,N=5\quad\quad$}},
legend style={at={(0.0,0.0)}, anchor=south west, draw=black}
]
\addlegendimage{no markers, red}
\addlegendimage{no markers, green!50.196078431372548!black}
\addlegendimage{no markers, blue}
\addlegendimage{no markers, red,dotted}
\addlegendimage{no markers, green!50.196078431372548!black,dotted}
\addlegendimage{no markers, blue,dotted}
\addlegendimage{no markers, red,dashed}
\addlegendimage{no markers, green!50.196078431372548!black,dashed}
\addlegendimage{no markers, blue,dashed}
\end{axis}
\end{tikzpicture}
	\caption{Solution at time $t = 0.14$ for the shock tube problem using knudsen number $\text{kn}=0.01$. The first 2 rows show the macroscopic qunatities $\rho_f,\Vvec_f,T_f,p_f$. In the thrid row we show $\Vvec$ and $T$ for different orders and different values of $c$ in \refer{eq:variational_formulation_V_and_T}.}
	\label{fig:shocktube_kn001}

\end{figure}

\begin{figure}
    \centering
    \begin{subfigure}[b]{0.49\textwidth}
		\resizebox{0.999\linewidth}{!}{\input{./images/shock_tube/compare_art_vs/rho_200_els_kn0.001_artvs_2.5e-05.tikz}}
    \end{subfigure}
    \hfill
    \begin{subfigure}[b]{0.49\textwidth}
		\resizebox{0.999\linewidth}{!}{\input{./images/shock_tube/compare_art_vs/pressure_200_els_kn0.001_artvs_2.5e-05.tikz}}
    \end{subfigure}

	\vspace{0.5cm}    
        
    \begin{subfigure}[b]{0.49\textwidth}
		\resizebox{0.999\linewidth}{!}{\input{./images/shock_tube/compare_art_vs/temperature_200_els_kn0.001_artvs_2.5e-05.tikz}}
    \end{subfigure}
    \hfill
    \begin{subfigure}[b]{0.49\textwidth}
		\resizebox{0.999\linewidth}{!}{\input{./images/shock_tube/compare_art_vs/v_200_els_kn0.001_artvs_2.5e-05.tikz}}
    \end{subfigure}

	\vspace{0.5cm}    
        
    \begin{subfigure}[b]{0.49\textwidth}
		\resizebox{0.999\linewidth}{!}{\input{./images/shock_tube/compare_art_vs/temperatureansatz_200_els_kn0.001_artvs_2.5e-05.tikz}}
    \end{subfigure}
    \hfill
    \begin{subfigure}[b]{0.49\textwidth}
		\resizebox{0.999\linewidth}{!}{\input{./images/shock_tube/compare_art_vs/vansatz_200_els_kn0.001_artvs_2.5e-05.tikz}}
    \end{subfigure}
   	\begin{tikzpicture}

\definecolor{color2}{rgb}{0.172549019607843,0.627450980392157,0.172549019607843}
\definecolor{color1}{rgb}{1,0.498039215686275,0.0549019607843137}
\definecolor{color4}{rgb}{0.580392156862745,0.403921568627451,0.741176470588235}
\definecolor{color0}{rgb}{0.12156862745098,0.466666666666667,0.705882352941177}
\definecolor{color3}{rgb}{0.83921568627451,0.152941176470588,0.156862745098039}

\begin{axis}
[
hide axis,
xmin=-2.995, xmax=62.895,
ymin=5.53988011611066e-19, ymax=0.216926593216958,
ymode=log,
tick align=outside,
tick pos=left,
x grid style={lightgray!92.026143790849673!black},
legend cell align={left},
y grid style={lightgray!92.026143790849673!black},
legend columns=3,
legend entries={{$c=64,N=3\quad\quad$},{$c=16,N=3\quad\quad$},{$c=4,N=3\quad\quad$},{$c=64,N=4\quad\quad$},{$c=16,N=4\quad\quad$},{$c=4,N=4\quad\quad$},{$c=64,N=5\quad\quad$},{$c=16,N=5\quad\quad$},{$c=4,N=5\quad\quad$}},
legend style={at={(0.0,0.0)}, anchor=south west, draw=black}
]
\addlegendimage{no markers, red}
\addlegendimage{no markers, green!50.196078431372548!black}
\addlegendimage{no markers, blue}
\addlegendimage{no markers, red,dotted}
\addlegendimage{no markers, green!50.196078431372548!black,dotted}
\addlegendimage{no markers, blue,dotted}
\addlegendimage{no markers, red,dashed}
\addlegendimage{no markers, green!50.196078431372548!black,dashed}
\addlegendimage{no markers, blue,dashed}
\end{axis}
\end{tikzpicture}
	\caption{Solution at time $t = 0.14$ for the shock tube problem using knudsen number $\text{kn}=0.001$. The first 2 rows show the macroscopic qunatities $\rho_f,\Vvec_f,T_f,p_f$. In the thrid row we show $\Vvec$ and $T$ for different orders and different values of $c$ in \refer{eq:variational_formulation_V_and_T}.}
	\label{fig:shocktube_kn0001}
\end{figure}

\begin{figure}
\centering
\resizebox{0.49\linewidth}{!}{
\begin{tikzpicture}
\begin{axis}[
hide axis,
scale only axis,
height = 0pt,
width = 0pt,
colormap name = ngsolve,
colorbar sampled,
colorbar horizontal,
point meta min = 0.125,
point meta max = 3.3,
colorbar style ={samples = 9,width = 0.55\linewidth, height = 0.3cm, xtick={0.125,1.7125,3.3}, xticklabel style = {style={font=\footnotesize}},
xticklabel pos = top
}]
\addplot [draw = none] coordinates { (0,0) };
\end{axis}
\end{tikzpicture}
}
\resizebox{0.49\linewidth}{!}{
\begin{tikzpicture}
\begin{axis}[
hide axis,
scale only axis,
height = 0pt,
width = 0pt,
colormap name = ngsolve,
colorbar sampled,
colorbar horizontal,
point meta min = 0.2,
point meta max = 4.2,
colorbar style ={samples = 9,width = 0.55\linewidth, height = 0.3cm, xtick={0.2,2.2,4.2}, xticklabel style = {style={font=\footnotesize}},
xticklabel pos = top
}]
\addplot [draw = none] coordinates { (0,0) };
\end{axis}
\end{tikzpicture}
}
	\vspace{0.25cm}

    \centering
    \begin{subfigure}[b]{0.98\textwidth}
		\resizebox{0.49\linewidth}{!}{\includegraphics{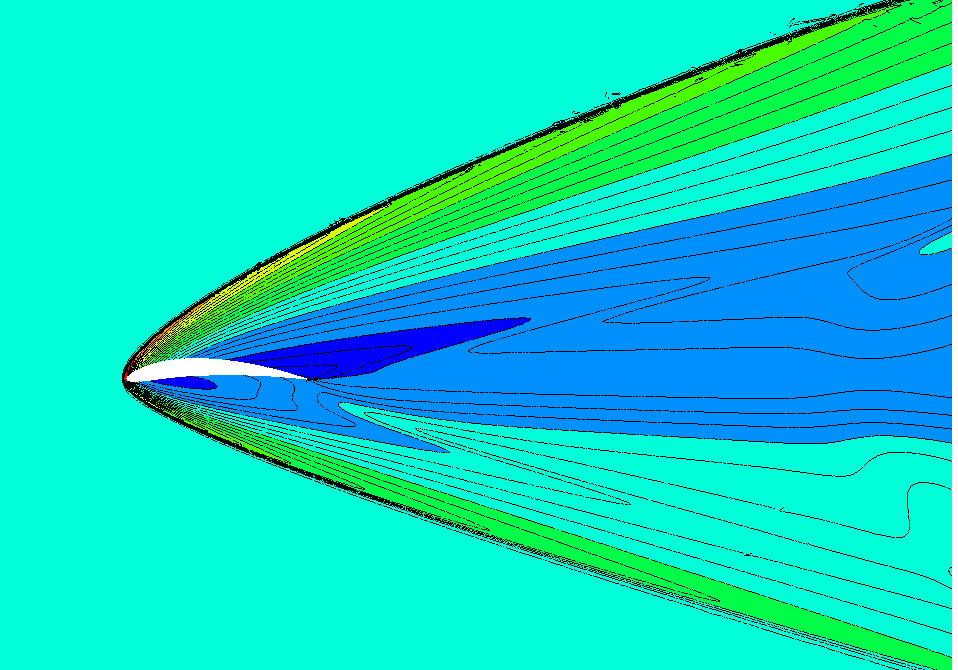}}
	\hfill
		\resizebox{0.49\linewidth}{!}{\includegraphics{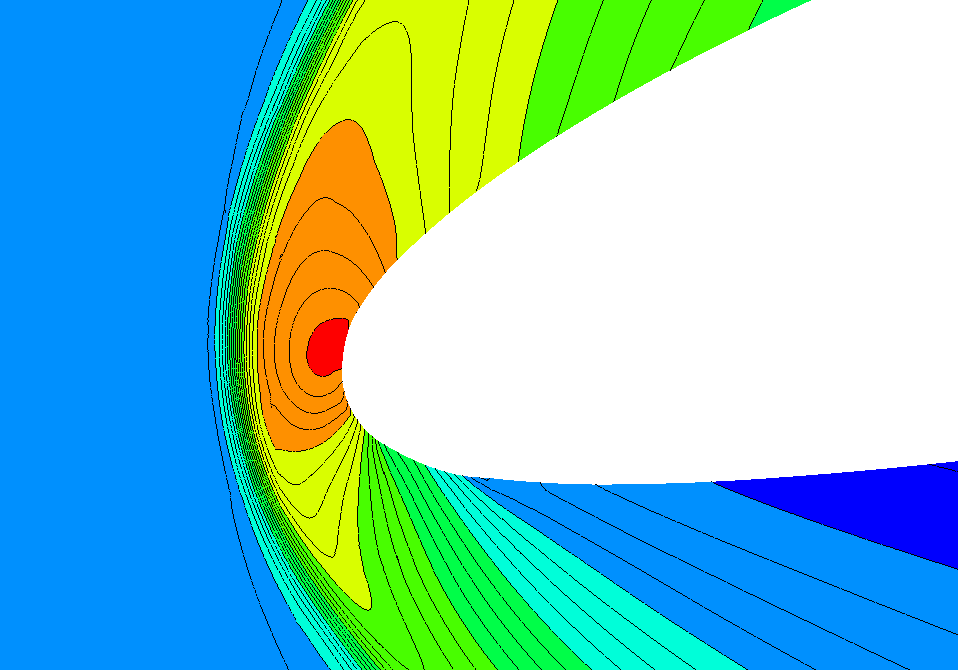}}
	\subcaption{Density $\rho_f$}
	\label{fig:rho_bow}
    \end{subfigure}
    
    \vspace{0.2cm}
\resizebox{0.49\linewidth}{!}{
\begin{tikzpicture}
\begin{axis}[
hide axis,
scale only axis,
height = 0pt,
width = 0pt,
colormap name = ngsolve,
colorbar sampled,
colorbar horizontal,
point meta min = 5,
point meta max = 10,
colorbar style ={samples = 9,width = 0.55\linewidth, height = 0.3cm, xtick={5,7.5,10}, xticklabel style = {style={font=\footnotesize}},
xticklabel pos = top
}]
\addplot [draw = none] coordinates { (0,0) };
\end{axis}
\end{tikzpicture}
}
\resizebox{0.49\linewidth}{!}{
\begin{tikzpicture}
\begin{axis}[
hide axis,
scale only axis,
height = 0pt,
width = 0pt,
colormap name = ngsolve,
colorbar sampled,
colorbar horizontal,
point meta min = 3.5,
point meta max = 28,
colorbar style ={samples = 9,width = 0.55\linewidth, height = 0.3cm, xtick={3.5,15.75,28}, xticklabel style = {style={font=\footnotesize}},
xticklabel pos = top
}]
\addplot [draw = none] coordinates { (0,0) };
\end{axis}
\end{tikzpicture}
}
	\vspace{0.25cm}

    \centering
    \begin{subfigure}[b]{0.98\textwidth}
		\resizebox{0.49\linewidth}{!}{\includegraphics{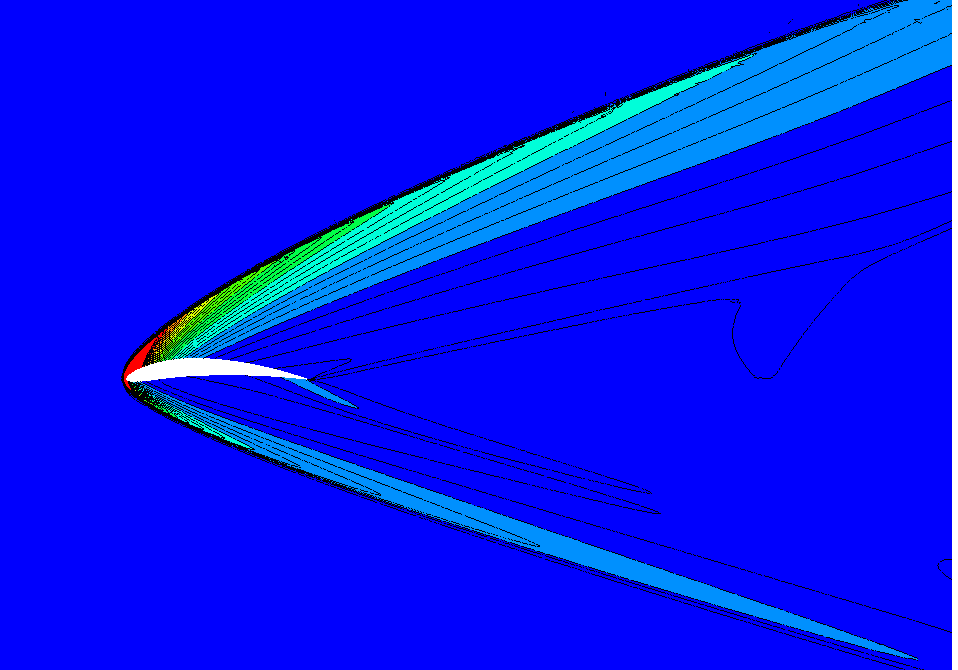}}
    \hfill    
		\resizebox{0.49\linewidth}{!}{\includegraphics{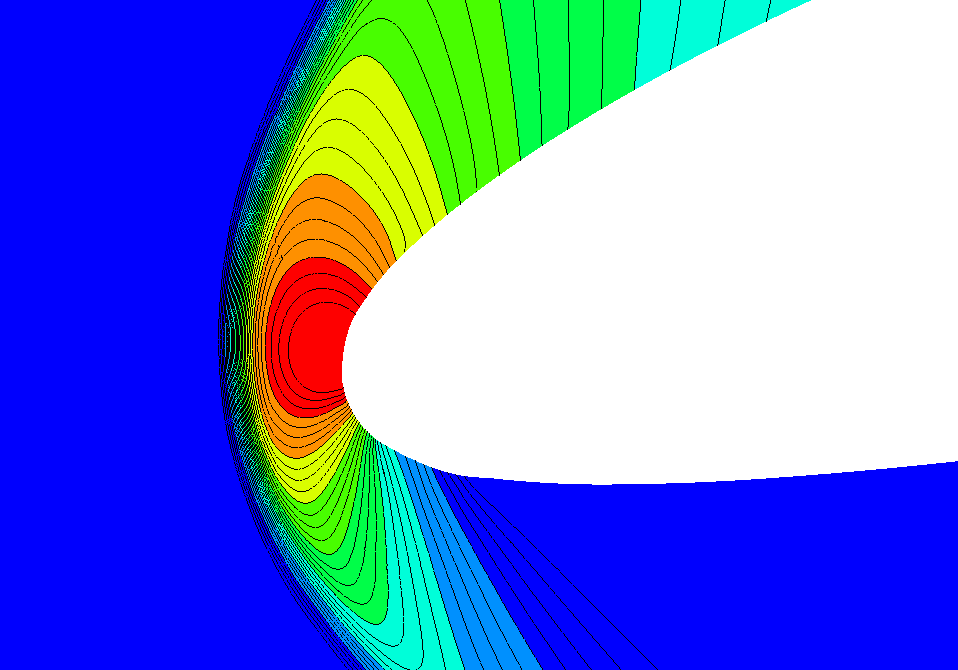}}
	\subcaption{presssure $p_f$}
	\label{fig:p_bow}
    \end{subfigure}
\end{figure}
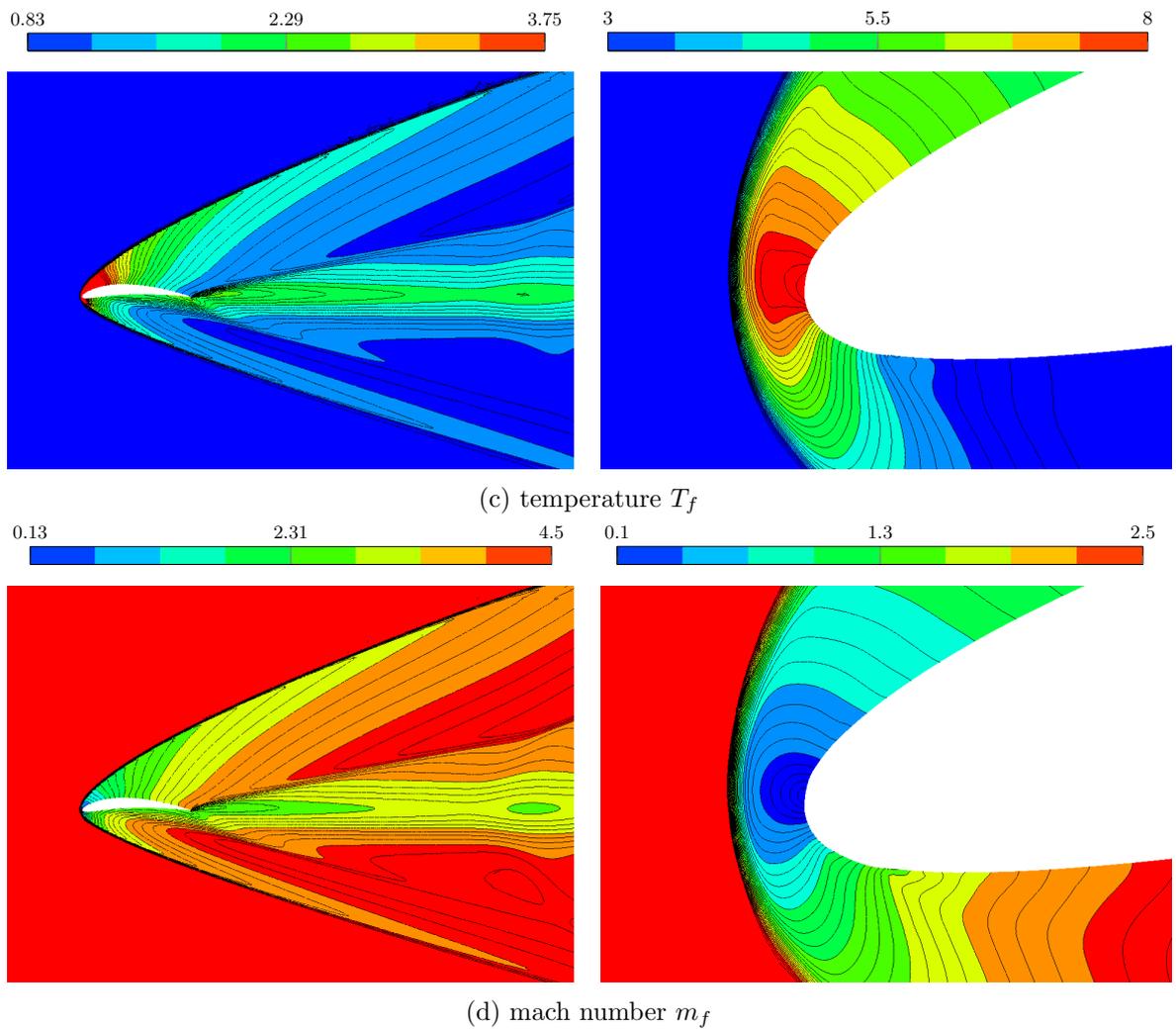
\begin{figure}    
\ContinuedFloat
\resizebox{0.49\linewidth}{!}{
\begin{tikzpicture}
\begin{axis}[
hide axis,
scale only axis,
height = 0pt,
width = 0pt,
colormap name = ngsolve,
colorbar sampled,
colorbar horizontal,
point meta min = 0.825,
point meta max = 3.75,
colorbar style ={samples = 9,width = 0.55\linewidth, height = 0.3cm, xtick={0.825,2.2875,3.75}, xticklabel style = {style={font=\footnotesize}},
xticklabel pos = top
}]
\addplot [draw = none] coordinates { (0,0) };
\end{axis}
\end{tikzpicture}
}
\resizebox{0.49\linewidth}{!}{
\begin{tikzpicture}
\begin{axis}[
hide axis,
scale only axis,
height = 0pt,
width = 0pt,
colormap name = ngsolve,
colorbar sampled,
colorbar horizontal,
point meta min = 3,
point meta max = 8,
colorbar style ={samples = 9,width = 0.55\linewidth, height = 0.3cm, xtick={3,5.5,8}, xticklabel style = {style={font=\footnotesize}},
xticklabel pos = top
}]
\addplot [draw = none] coordinates { (0,0) };
\end{axis}
\end{tikzpicture}
}
	\vspace{0.25cm}
	
    \centering
    \begin{subfigure}[b]{0.98\textwidth}
		\resizebox{0.49\linewidth}{!}{\includegraphics{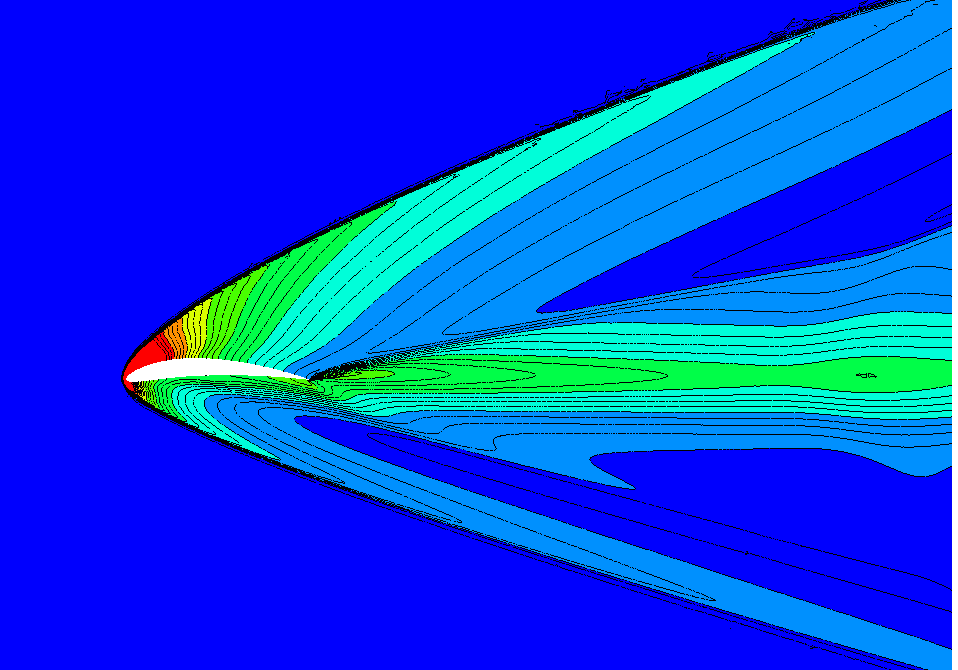}}
    \hfill    
		\resizebox{0.49\linewidth}{!}{\includegraphics{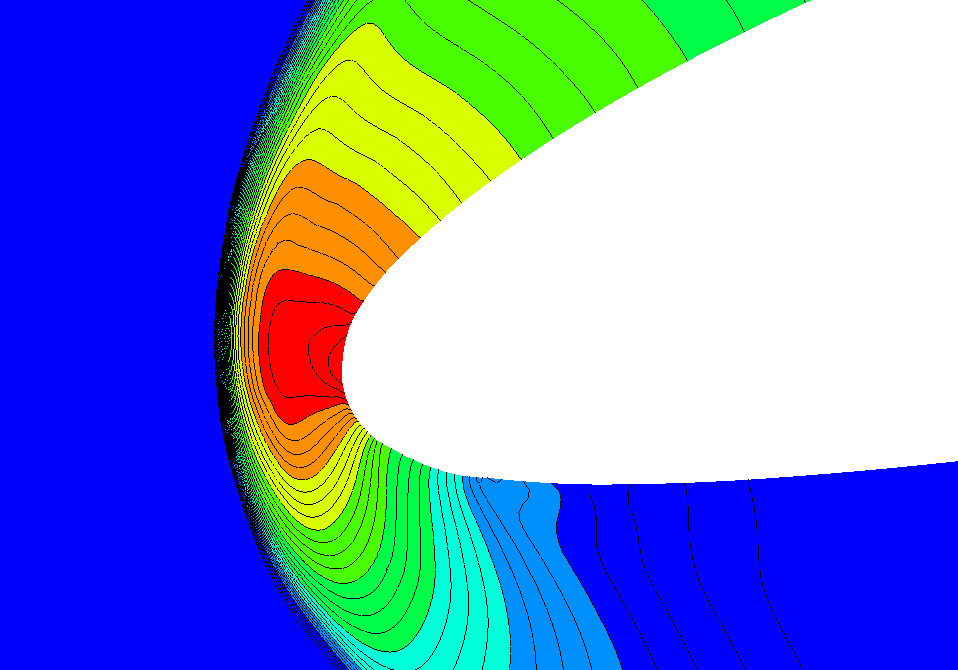}}
	\subcaption{temperature $T_f$}
	\label{fig:temp_bow}
    \end{subfigure}
\resizebox{0.49\linewidth}{!}{
\begin{tikzpicture}
\begin{axis}[
hide axis,
scale only axis,
height = 0pt,
width = 0pt,
colormap name = ngsolve,
colorbar sampled,
colorbar horizontal,
point meta min = 0.125,
point meta max = 4.5,
colorbar style ={samples = 9,width = 0.55\linewidth, height = 0.3cm, xtick={0.125,2.3125,4.5}, xticklabel style = {style={font=\footnotesize}},
xticklabel pos = top
}]
\addplot [draw = none] coordinates { (0,0) };
\end{axis}
\end{tikzpicture}
}
\resizebox{0.49\linewidth}{!}{
\begin{tikzpicture}
\begin{axis}[
hide axis,
scale only axis,
height = 0pt,
width = 0pt,
colormap name = ngsolve,
colorbar sampled,
colorbar horizontal,
point meta min = 0.1,
point meta max = 2.5,
colorbar style ={samples = 9,width = 0.55\linewidth, height = 0.3cm, xtick={0.1,1.3,2.5}, xticklabel style = {style={font=\footnotesize}},
xticklabel pos = top
}]
\addplot [draw = none] coordinates { (0,0) };
\end{axis}
\end{tikzpicture}
}
    \vspace{0.25cm}
    
    \centering
    \begin{subfigure}[b]{0.98\textwidth}
		\resizebox{0.49\linewidth}{!}{\includegraphics{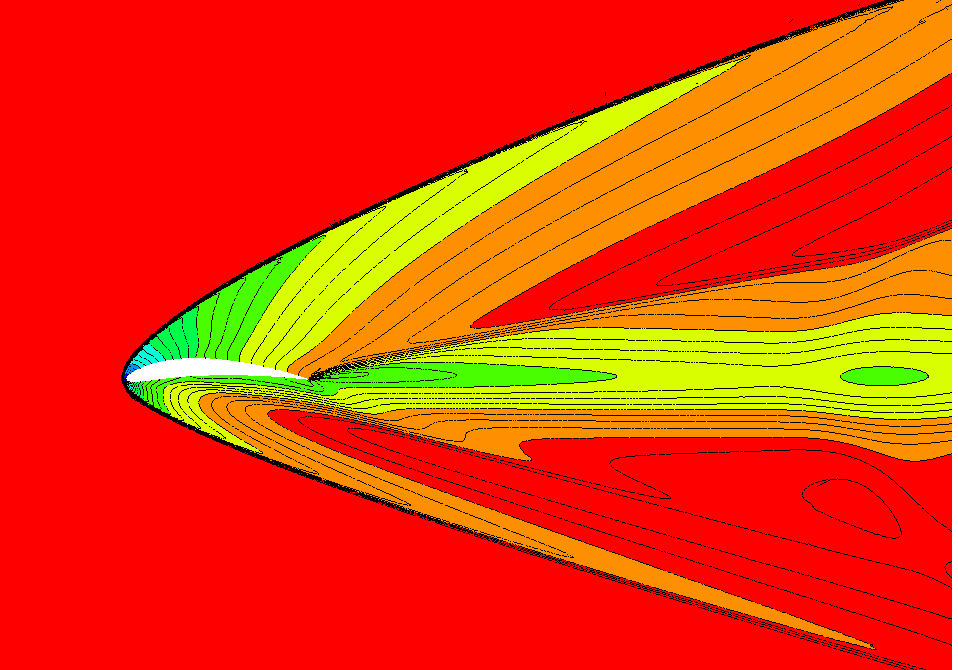}}
    \hfill    
		\resizebox{0.49\linewidth}{!}{\includegraphics{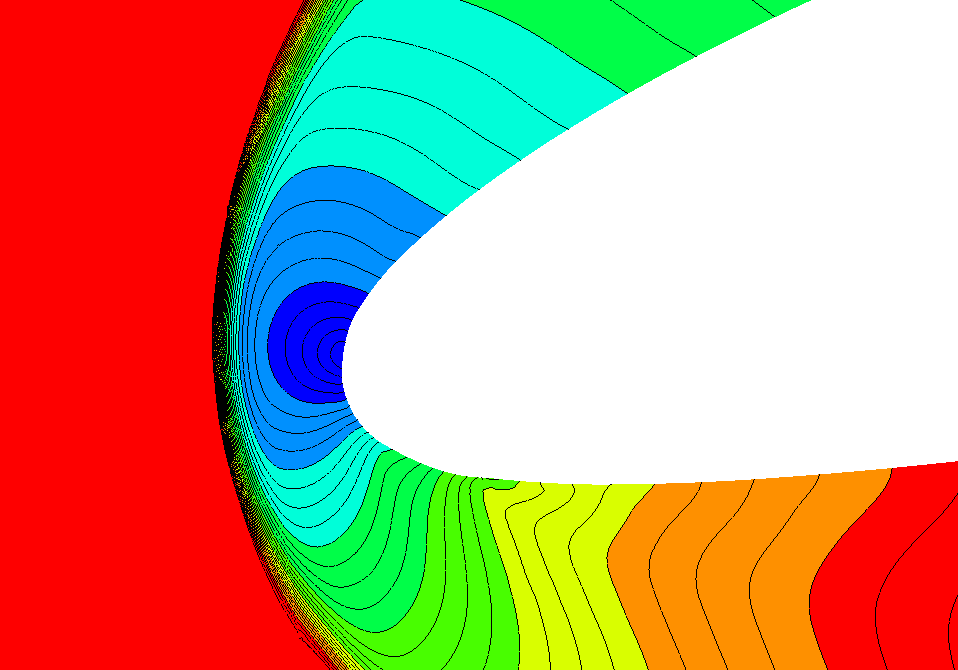}}
	\subcaption{mach number $m_f$}
	\label{fig:Mach_bow}		
    \end{subfigure}

	\caption{Macroscopic fields for a Mach 4.5 flow around a NACA 7410 air foil, example \ref{sec:numerical_results2}. Calculations were done with order 6 polynomials for $\xvec$ and order 3 for $\vvec$. Time stepping with forward euler scheme, using $\tau = 7.8125\text{e-}6$. The smoothing parameter in  \refer{eq:variational_formulation_V_and_T} was $c=18$.}
	\label{fig:NACA7410}
\end{figure}

\bibliographystyle{abbrv}
\FloatBarrier

\end{document}